\documentclass[11pt,a4paper]{article}
\usepackage{amsfonts,amsmath,amssymb,accents,amsthm,bbm}
\usepackage{verbatim,bm}
\usepackage{hyperref}
\usepackage{tikz}
\usetikzlibrary{calc,arrows,positioning,fit,petri,external}

\newcommand{\dis}{\displaystyle}

\newcommand{\noi}{\noindent}
\newcommand{\halmos}{\rule{1ex}{1.4ex}}
\newcommand{\QED}{\nopagebreak{\hspace*{\fill}$\halmos$\medskip}}

\newcommand{\med}{\medskip}
\newcommand{\quand}{\quad\mbox{and}\quad}

\newtheoremstyle{mythm}
{}
{}
{\itshape}
{}
{\bfseries}
{}
{.5em}
{#1 #2 \thmnote{(#3)}}

\newtheoremstyle{mydef}
{}
{}
{\upshape}
{}
{\bfseries}
{}
{.5em}
{#1 #2 \thmnote{(#3)}}
  
\theoremstyle{mythm}

\newtheorem{theorem}{Theorem}
\newtheorem{proposition}[theorem]{Proposition}
\newtheorem{lemma}[theorem]{Lemma}
\newtheorem{conjecture}[theorem]{Conjecture}
\theoremstyle{mydef}
\newtheorem{remark}[theorem]{Remark}
\newtheorem{defi}[theorem]{Definition}
\newcommand{\bt}{\begin{theorem}}
\newcommand{\et}{\end{theorem}}
\newcommand{\bl}{\begin{lemma}}
\newcommand{\el}{\end{lemma}}
\newcommand{\bp}{\begin{proposition}}
\newcommand{\ep}{\end{proposition}}
\newcommand{\br}{\begin{remark}}
\newcommand{\er}{\end{remark}}
\newcommand{\bcon}{\begin{conjecture}}
\newcommand{\econ}{\end{conjecture}}

\newenvironment{Proof}[1][]{\noi\textbf{Proof #1}}{\QED}
\newcommand{\bpro}{\begin{Proof}}
\newcommand{\epro}{\end{Proof}}

\newcommand{\be}{\begin{equation}}
\newcommand{\ee}{\end{equation}}
\newcommand{\ba}{\begin{array}}
\newcommand{\ea}{\end{array}}
\newcommand{\bac}{\begin{array}{r@{\,}c@{\,}l}}
\newcommand{\bc}{\be\begin{array}{r@{\,}c@{\,}l}}
\newcommand{\ec}{\end{array}\ee}
\newcommand{\bcase}{\left\{\ba{l@{\quad\mbox{if }}l}}
\newcommand{\ecase}{\ea\right.}

\expandafter\mathchardef\expandafter\varphi\number\expandafter\phi\expandafter\relax
\expandafter\mathchardef\expandafter\phi\number\varphi

\newcommand{\de}{\delta}
\newcommand{\De}{\Delta}
\newcommand{\eps}{\varepsilon}
\newcommand{\La}{\Lambda}
\newcommand{\sig}{\sigma}
\newcommand{\phh}{\varphi}
\newcommand{\Phh}{\Phi}
\newcommand{\phhb}{{\bm{\phh}}}
\newcommand{\psib}{{\bm{\psi}}}

\newcommand{\Ei}{{\cal E}}
\newcommand{\Fi}{{\cal F}}
\newcommand{\Gi}{{\cal G}}
\newcommand{\Hi}{{\cal H}}
\newcommand{\Oi}{{\cal O}}
\newcommand{\Ti}{{\cal T}}

\newcommand{\Vi}{{\cal V}}
\newcommand{\Wi}{{\cal W}}
\newcommand{\Zi}{{\cal Z}}
\newcommand{\N}{{\mathbb N}}
\renewcommand{\P}{{\mathbb P}}
\newcommand{\R}{{\mathbb R}}
\newcommand{\Z}{{\mathbb Z}}

\newcommand{\pre}{\preceq}
\newcommand{\desd}{\ensuremath{\Leftrightarrow}}

\newcommand{\sub}{\subset}
\newcommand{\beh}{\backslash}

\newcommand{\asto}[1]{\underset{{#1}\to\infty}{\longrightarrow}}
\newcommand{\Asto}[1]{\underset{{#1}\to\infty}{\Longrightarrow}}
\newcommand{\ti}{\tilde}
\newcommand{\ov}{\overline}
\newcommand{\un}{\underline}
\newcommand{\pa}{\partial}
\newcommand{\ffrac}[2]{{\textstyle\frac{{#1}}{{#2}}}}
\newcommand{\ha}{\ffrac{1}{2}}

\newcommand{\rbf}{\mathbf{r}}

\newcommand{\lis}{1\leq s\leq\sig}

\begin{document}

\makeatletter\@addtoreset{equation}{section}
\makeatother\def\theequation{\thesection.\arabic{equation}} 

\renewcommand{\labelenumi}{{\rm (\roman{enumi})}}
\renewcommand{\theenumi}{\roman{enumi}}

\title{Peierls bounds from Toom contours}

\author{Jan M.~Swart\footnote{The Czech Academy of Sciences, Institute of Information Theory and Automation, Pod vod\'arenskou v\v{e}\v{z}\'i 4, 18200 Praha 8, Czech republic. swart@utia.cas.cz},
R\'eka Szab\'o\footnote{Bernoulli Institute, University of Groningen, Nijenborgh 9, 9747 AG Groningen, The Netherlands. r.szabo@rug.nl},
and Cristina Toninelli\footnote{CEREMADE, CNRS, Universit\'e Paris-Dauphine, PSL University, Place du Mar\'echal de Lattre de Tassigny, 75016 Paris, France and DMA, \'Ecole normale sup\'erieure, PSL University,  45 rue d’Ulm, F-75230 Cedex 5 Paris, France. toninelli@ceremade.dauphine.fr}}

\date{\today}

\maketitle

\begin{abstract}\noi
For deterministic monotone cellular automata on the $d$-dimensional integer lattice, Toom has given necessary and sufficient conditions for the all-one fixed point to be stable against small random perturbations. The proof of sufficiency is based on an intricate Peierls argument. We present a simplified version of this Peierls argument. Our main motivation is the open problem of determining stability of monotone cellular automata with intrinsic randomness, in which for the unperturbed evolution the local update rules at different space-time points are chosen in an i.i.d.\ fashion according to some fixed law. We apply Toom's Peierls argument to prove stability of a class of cellular automata with intrinsic randomness and also derive lower bounds on the critical parameter for some deterministic cellular automata.
\end{abstract}
\vspace{.5cm}

\noi

\newpage

{\setlength{\parskip}{-2pt}\tableofcontents}

\newpage

\part{Results}

\section{Introduction and main results}\label{S:intmain}

\subsection{Introduction}

A cellular automaton is a discrete-time dynamical process on the set of functions from the $d$-dimensional integer lattice to some finite local state space. In each time step, the state of each point in the lattice is updated according to an update rule that uses information from finitely many near-by lattice points to determine the new state of a given lattice point. Let us assume that the update rules in all lattice points are the same up to translation, that the local state space is $\{0,1\}$, and that the update rule is monotone in the sense that it preserves the order (maps configurations with more 1's into configurations with more 1's). Toom \cite{Too80} studied random perturbations of such monotone cellular automata. In his celebrated stability theorem \cite[Thm~5]{Too80}, he gave necessary and sufficient conditions, in terms of the update rule, for the all-one fixed point of a monotone cellular automaton to be stable against small random perturbations.

In this paper, we are interested in monotone cellular automata for which the unperturbed dynamics are already random, in the sense that the update rules that are used in different space-time points are chosen from some finite set of possible update rules in an i.i.d.\ fashion. Also for such monotone cellular automata with ``intrinsic randomness'', we ask whether the all-one fixed point is stable against small random perturbations. Proving stability for monotone cellular automata with intrinsic randomness is much harder than in the deterministic case. In Theorem~\ref{T:detspeed}, we give sufficient conditions for stability, that are, however, far from necessary. Nevertheless, Theorem~\ref{T:detspeed} contains Toom's stability theorem as a special case, demonstrating the fact that the case with intrinsic randomness is harder. In a companion paper \cite{SSTcont}, building further on the theory developed in the present paper, we will be able to improve on Theorem~\ref{T:detspeed} while still not completely solving the stability problem for monotone cellular automata with intrinsic randomness.

The main contribution of the present paper lies not so much in Theorem~\ref{T:detspeed} itself as in our methods to prove it. We base ourselves on the original Peierls argument that Toom used to prove his stability theorem, but spend a lot of time reformulating and reorganising the argument, which culminates in Theorems \ref{T:contour}, \ref{T:cycle}, and \ref{T:ToomPei}. Theorem~\ref{T:contour}, in particular, makes explicit and significantly reformulates a statement that has so far remained hidden in the proofs of \cite{Too80}, and that really is the cornerstone of his Peierls argument. By making this explicit, we hope to facilitate the use of this technique in future. This has to some degree already happened, as these techniques have already been applied not only in our companion paper \cite{SSTcont}, but also in \cite{HS22} which is about bootstrap percolation and in \cite{CSS24} which is about random games.

\subsection{Set-up and background}\label{S:intro}

Let $\{0,1\}^{\Z^d}$ denote the set of configurations $x=(x(i))_{i\in\Z^d}$ of zeros and ones on the $d$-dimensional integer lattice $\Z^d$. By definition, a map $\phi:\{0,1\}^{\Z^d}\to\{0,1\}$ is \emph{local} if $\phi$ depends only on finitely many coordinates, i.e., there exists a finite set $\De\sub\Z^d$ and a function $\phi':\{0,1\}^\De\to\{0,1\}$ such that $\phi\big((x(i))_{i\in\Z^d}\big)=\phi'\big((x(i))_{i\in\De}\big)$ for each $x\in\{0,1\}^{\Z^d}$. We let $\De(\phi)$ denote the smallest set with this property, which may be empty: in this case $\phi$ is either constantly zero or one. We denote the constant functions by
\be\label{consmap}
\phi^0(x):=0\quand\phi^1(x):=1\qquad(x\in\{0,1\}^{\Z^d}).
\ee
A local map $\phi$ is \emph{monotone} if $x\leq y$ (coordinatewise) implies $\phi(x)\leq\phi(y)$. Let $\{\phi_0,\ldots,\phi_m\}$ be a set of monotone local maps $\phi_k:\{0,1\}^{\Z^d}\to\{0,1\}$, of which $\phi_0=\phi^0$ is the map that is constantly zero and $\phi_1,\ldots,\phi_m$ are not constant. Let $\rbf=\big(\rbf(1),\ldots,\rbf(m)\big)$ be a probability distribution on $\{1,\ldots,m\}$. We will be interested in i.i.d.\ collections of random variables $\Phi^{p,\rbf}=\Phi^p=(\Phi^p_{i,t})_{(i,t)\in\Z^{d+1}}$ with values in $\{\phi_0,\ldots,\phi_m\}$ such that
\be\label{Phip}
\P\big[\Phi^p_{i,t}=\phi_k\big]=\left\{\ba{ll}
\dis p\quad&\mbox{if }k=0,\\[5pt]
\dis (1-p)\rbf(k)\quad&\mbox{if }1\leq k\leq m,
\ea\right.
\ee
where $p\in[0,1]$ is a parameter. We call $\Phi^p$ a \emph{monotone cellular automaton}. We will be interested in the case that $p$ is small but positive. We think of $\Phi^p$ as a small perturbation of $\Phi^0$. In the special case that $m=1$, we say that $\Phi^0$ is a \emph{deterministic} monotone cellular automaton. If $m\geq 2$ and $\rbf(k)<1$ for all $k$, then we say that $\Phi^0$ has \emph{intrinsic randomness}.

If $X^p_0$ is a random variable with values in $\{0,1\}^{\Z^d}$, independent of $(\Phi^p_{i,t})_{i\in\Z^d,\ t\in\mathbb Z_+}$, then setting
\be\label{Markov}
X^p_t(i):=\Phi^p_{i,t}\big((X^p_{t-1}(i+j))_{j\in\Z^d}\big)\qquad(i\in\Z^d,\ t>0)
\ee
defines a Markov chain $(X^p_t)_{t\geq 0}$ with state space $\{0,1\}^{\Z^d}$. In words, the new state of $i$ at time $t$ is obtained by applying the random map $\Phi^p_{i,t}$ to the configuration at time $t-1$, shifted so that $i$ is located at the origin. Let $\P^x$ denote the law of this Markov chain started in a given initial state $X^p_0=x$ and let $\un 0$ and $\un 1$ denote the configurations in $\{0,1\}^{\Z^d}$ that are constantly zero or one, respectively. It is well-known (compare \cite[Thm~III.2.3]{Lig85}) that
\be\label{toupper}
\P^{\un 0}\big[X^p_t\in\,\cdot\,\big]\Asto{t}\un\nu_p
\quand
\P^{\un 1}\big[X^p_t\in\,\cdot\,\big]\Asto{t}\ov\nu_p
\ee
where $\Rightarrow$ denotes weak convergence of probability measures on $\{0,1\}^{\Z^d}$, equipped with the product topology, and $\un\nu_p$ and $\ov\nu_p$ are invariant laws of the Markov chain defined in (\ref{Markov}), that are called the \emph{lower} and \emph{upper} invariant laws, respectively. Let
\be\label{ovrho}
\ov\rho_\rbf(p)=\ov\rho(p):=\lim_{t\to\infty}\P^{\un 1}\big[X^p_t(i)=1\big]\qquad(p\in[0,1],\ i\in\Z^d)
\ee
denote the density of the upper invariant law, which by translation invariance does not depend on $i\in\Z^d$. Trivially, $\ov\rho(0)=1$ \footnote{Indeed from the fact that  $\phi_k$ are monotone and non constant for $k\in\{1,\dots,m\}$,  it follows that $\phi_k(\un 1)=1$.}.
We say that the monotone cellular automaton $\Phi^0$ defined by the monotone local maps $\phi_1,\ldots,\phi_m$ and the probability distribution $\rbf$ is \emph{stable} if
\be\label{stable}
\lim_{p\to 0}\ov\rho(p)=1,
\ee
and \emph{completely unstable} if $\ov\rho(p)=0$ for all $p>0$. 
A simple coupling argument shows that $p\mapsto\ov\rho(p)$ is non-increasing, so
if we let
\be\label{pc}
p_{\rm c}:=\sup\big\{p\in[0,1]:\ov\rho(p)>0\big\}
\ee
it holds that $\ov\rho(p)>0$ for all $p<p_{\rm c}$ and $\ov\rho(p)=0$ for all $p>p_{\rm c}$. In particular, complete instability corresponds to $p_c=0$.

For deterministic monotone cellular automata, Toom \cite{Too80} has completely solved the problem of determining whether a given cellular automaton is stable or not. To state his result we first need to define \emph{eroders}. For each local map $\phi:\{0,1\}^{\Z^d}\to\{0,1\}$, we let $\Psi_\phi:\{0,1\}^{\Z^d}\to\{0,1\}^{\Z^d}$ be defined as
\be\label{Psiphi}
\Psi_\phi(x)(i):=\phi\big((x(i+j))_{j\in\Z^d}\big)\qquad\big(x\in\{0,1\}^{\Z^d}\big),
\ee
i.e., $\Psi_\phi$ describes one step of the time evolution of the deterministic cellular automaton defined by $\phi$. 

\begin{defi}[Eroders]\label{def:eroders}
We say that a local map $\phi$ is an \emph{eroder} if for each configuration $x\in\{0,1\}^{\Z^d}$ that contains only finitely many zeros, there is a $t\in\N$ such that $\Psi^t_\phi(x)=\un 1$, where $\Psi^t_\phi$ denotes the $t$-th iterate of the map $\Psi_\phi$.
\end{defi}

We quote the following result from \cite[Thm~5]{Too80}.\footnote{In the case where $\phi$ is not an eroder, Theorem~5 in \cite{Too80} only states that the monotone cellular automaton is not stable, but the proof actually implies that it is completely unstable. We will give our own proof in Subsection~\ref{S:instab} below.}

\begin{theorem}[Toom's stability theorem]
The\label{T:Toom} deterministic monotone cellular automaton $\Phi^0$ defined by a monotone local nonconstant map $\phi$ is stable if $\phi$ is an eroder and completely unstable if $\phi$ is not an eroder.
\end{theorem}

For general local maps that need not be monotone, it is known that there exists no algorithm to decide whether a given map is an eroder, even in one dimension \cite{Pet87}. By contrast, for monotone local maps, there exists a simple criterion to check whether a given map is an eroder. To state this criterion we need the notion of minimal one-sets.
A \emph{one-set} of a monotone local map $\phi:\{0,1\}^{\Z^d}\to\{0,1\}$ is a finite set $A\sub\Z^d$ such that $\phi(1_A)=1$, where $1_A$ denotes the indicator function of $A$. A \emph{minimal one-set} is a one-set that does not contain other one-sets as a proper subset. Each monotone local map $\phi:\{0,1\}^{\Z^d}\to\{0,1\}$ can be written as
\be\label{Aphi}
\phi(x)=\bigvee_{A\in\Oi(\phi)}\bigwedge_{i\in A}x(i)\qquad\big(x\in\{0,1\}^{\Z^d}\big),
\ee
where $\Oi(\phi)$ is the set of minimal one-sets of $\phi$, and $\vee$ and $\wedge$ denote the maximum and minimum operations, respectively. In (\ref{Aphi}), we use the convention that the supremum (resp.\ infimum) over an empty set is 0 (resp.\ 1). In line with this, $\Oi(\phi^0)=\emptyset$ and $\Oi(\phi^1)=\{\emptyset\}$ (note the difference!). We let ${\rm Conv}(A)$ denote the convex hull of a set $A$, viewed as a subset of $\R^d$. Then \cite[Thm~6]{Too80}, with simplifications due to \cite[Thm~1]{Pon13}, says the following.

\bp[Erosion criterion]\label{P:ero}
A monotone local map $\phi\neq\phi^0$ is an eroder if and only if
\be\label{erosion}
\bigcap_{A\in\Oi(\phi)}{\rm Conv}(A)=\emptyset.
\ee
\ep

See also Lemma~\ref{L:erode} which gives a related alternative erosion criterion due to \cite[Lemma~12]{Pon13}.

\begin{remark}\label{rem:Helly}
Helly's theorem \cite[Corollary~21.3.2]{Roc70} guarantees that if (\ref{erosion}) holds, then there exists a subset $\Oi'\sub\Oi(\phi)$ of cardinality at most $d+1$ such that $\bigcap_{A\in\Oi'}{\rm Conv}(A)=\emptyset$.
\end{remark}

For concreteness, let us look at some examples of maps in two dimensions. We set
\bc\label{phiNEC}
\dis\phi^{\rm NEC}(x)&:=&\dis{\tt round}\big([x(0,0)+x(0,1)+x(1,0)]/3\big),\\[5pt]
\dis\phi^{\rm NN}(x)&:=&\dis{\tt round}\big([x(0,0)+x(0,1)+x(1,0)+x(0,-1)+x(-1,0)]/5\big),\\[5pt]
\dis\phi^{\rm coop}(x)&:=&\dis x(0,0)\vee\big(x(0,1)\wedge x(1,0)\big),
\ec
where ${\tt round}$ denotes the function that rounds off a real number to the nearest integer. The function $\phi^{\rm NEC}$ is known as \emph{North-East-Centre voting} or \emph{NEC voting}, for short, and also as \emph{Toom's rule}. In analogy with $\phi^{\rm NEC}$, we also define maps $\phi^{\rm NWC},\phi^{\rm SWC},\phi^{\rm SEC}$ that describe North-West-Centre voting, South-West-Centre voting, and South-East-Centre voting, respectively, defined in the obvious way. We will call the map $\phi^{\rm NN}$ from (\ref{phiNEC}) \emph{Nearest Neighbour voting} or \emph{NN voting}, for short. Another name found in the literature is the \emph{symmetric majority rule}. We call $\phi^{\rm coop}$ the \emph{cooperative branching rule}. It is also known as the \emph{sexual reproduction rule} because of the interpretation that when $\phi^{\rm coop}$ is applied at a site $(i_1,i_2)$, two parents at $(i_1+1,i_2)$ and $(i_1,i_2+1)$ produce offspring at $(i_1,i_2)$, provided the parents' sites are both occupied and $(i_1,i_2)$ is vacant.
Using Proposition~\ref{P:ero} one can easily check that $\phi^{\rm NEC}$ and $\phi^{\rm coop}$ are eroders, but $\phi^{\rm NN}$ is not. Indeed, we have
\bc\label{ANEC}
\mathcal A(\dis\phi^{\rm NEC})&:=& \left\{\{(0,0), (1,0)\},\{(0,0), (0,1)\},\{(0,1),(1,0)\} \right\}, \\[5pt]
\mathcal A(\phi^{\rm coop})&:=&  \left\{\{(0,0)\},\{(0,1),(1,0)\} \right\} ,
\ec
and both sets satisfy condition \eqref{erosion}. On the other hand, $\mathcal A(\phi^{\rm NN})$ is the set of all subsets of cardinality 3 of $\{(0,0),(0,1),(1,0),(0,-1),(-1,0)\}$. Therefore, each $A\in\mathcal A(\phi^{\rm NN})$ contains the origin in its convex hull, and the erosion condition \eqref{erosion} is not satisfied. In fact, it is not hard to find configurations containing only finitely many zeros which cannot disappear under iterated applications of the map $\Psi_{\dis\phi^{\rm NN}}$, for example  the configuration that is zero on $(0,0),(0,1),(1,0),(1,1)$ and one everywhere else.

\begin{remark}
Toom's\label{R:multistep} stability theorem is stated in a slightly greater generality. The deterministic monotone cellular automata considered in \cite{Too80} are defined by monotone local maps $\phi$ that can ``look back'' more than one time step, in the sense that the set $\De(\phi)$ defined above (\ref{consmap}) is a finite subset of $\Z^d\times\Z_-$. In this case, $\phi$ is an eroder if and only if
\be
\bigcap_{A\in\Oi(\phi)}\bigcup_{\alpha>0}\{\alpha\cdot(i,t): (i, t)\in{\rm Conv}(A)\}=\emptyset,
\ee
that is, no ray in $\R^{d+1}$ that starts from the origin intersects all the convex hulls of the minimal one-sets. Note that in our setting this condition is equivalent to \eqref{erosion}.
\end{remark}

Toom's Theorem~\ref{T:Toom} settles the stability issue for deterministic monotone cellular automata. The next natural step is to study stability for monotone cellular automata with intrinsic randomness. One might think that stability should hold at least in the case when $\phi_1,\ldots,\phi_m$ are all eroders, but this is not true. For example, there are good reasons to believe that the monotone cellular automaton that applies the maps $\phi^{\rm NEC},\phi^{\rm NWC},\phi^{\rm SWC},\phi^{\rm SEC}$ each with probability $1/4$ is unstable, in spite of the fact that each of these maps individually is an eroder, see Conjecture~\ref{C:symNEC} below.

To see a further example of the difficulties of cellular automata with intrinsic randomness, consider the \emph{identity map}, defined as
\be\label{phiid}
\phi^{\rm id}(x):=x(0)\qquad\big(x\in\{0,1\}^{\Z^d}\big).
\ee
In terms of the associated Markov chain (\ref{Markov}), applying the identity map in a given space-time point has the effect that the local state at a site does not change. One might think that if $\phi$ is an eroder, then a cellular automaton that applies the maps $\phi$ and $\phi^{\rm id}$ each with positive probability must be stable, but again this turns out to be wrong. Gray \cite[Examples~18.3.5 and 18.3.6]{Gra99} has given convincing arguments that show that the addition of the identity map can make eroders unstable and conversely, make non-eroders stable. Being able to include the identity map is important for understanding continuous-time interacting particle systems. We can think of such systems as limits of discrete-time cellular automata where time is measured in steps of some small size $\de$ and all maps except $\phi^{\rm id}$ are applied with a probability of order~$\de$.

The most difficult part of Theorem~\ref{T:Toom} is the statement that $\Phi^0$ is stable if $\phi$ is an eroder. To prove this, Toom used an intricate Peierls argument. It is fair to say that Toom's original paper \cite{Too80} is quite hard to read. Indeed, several subsequent papers have been devoted to simplifying his arguments and others have re-proved his result from scratch for some specific model at hand to avoid relying on this complex proof \cite{LMS90,Gac95,Pre07,Pon13,Gac21} (see Subsection~\ref{S:discussion}).

In this paper, we reformulate and simplify Toom's Peierls argument. Our main motivation is the problem of extending Toom's stability theorem to monotone cellular automata with intrinsic randomness. As a first step in this direction, we will prove a stability result in Theorem~\ref{T:detspeed} below, which however excludes many interesting cases such as cellular automata that apply the identity map with a positive probability. This is not due to a fundamental limitation of Toom's Peierls argument, but to go beyond Theorem~\ref{T:detspeed} one needs more advanced methods to estimate the Peierls sum. In order not to overload the present paper, we have delegated these methods to a companion paper \cite{SSTcont} where further stability results for cellular automata with intrinsic randomness will be proved.

As a further result of our reformulation of Toom's Peierls argument, we will derive explicit lower bounds on the critical noise parameter $p_{\rm c}$ from (\ref{pc}) for some deterministic cellular automata. Although these bounds are often several orders of magnitude from the conjectured true values, they are nevertheless the sharpest rigorous bounds available. For a subclass of cellular automata, we show that it is possible to derive significantly better bounds by using Toom cycles, which are Toom contours with additional useful properties.

Toom's Peierls argument was invented to study stability of monotone cellular automata with respect to noise that is i.i.d.\ in space and time. It has recently been discovered that it can also be used to prove stability with respect to noise that is applied only to the initial state \cite{HS22,CSS24}. This has applications in bootstrap percolation, which we will briefly discuss in Appendix~\ref{A:boot} below.

\subsubsection*{Outline}

In the remainder of Section~\ref{S:intmain} we discuss applications of Toom's Peierls argument. Stability of monotone cellular automata with intrinsic randomness is discussed in Subsection~\ref{S:stab}, explicit lower bounds on the critical noise parameter are presented in Subsection~\ref{S:pc}. In Subsection~\ref{S:discussion} we discuss earlier work on the topic and state some open problems.

In Section~\ref{S:ToPe} we present our reformulation of Toom's Peierls argument. We work in a more general setting than in Section~\ref{S:intmain}, which also allows for cellular automata that can look back more than one time step as in Remark~\ref{R:multistep} and cellular automata on other lattices than $\Z^d$, such as trees. We show that each monotone cellular automaton has a maximal trajectory and that the density $\ov\rho(p)$ of the upper invariant law is equal to the probability that this maximal trajectory has a one at the origin. We moreover introduce objects we call \emph{Toom contours} that are directed graphs with different types of edges that are designed to make use of the characterisation of eroders in terms of edge speeds and polar functions that will be discussed in Subsection~\ref{S:stab} below.

The main results of Section~\ref{S:ToPe} and indeed of the whole paper are Theorems \ref{T:contour}, \ref{T:cycle}, and \ref{T:ToomPei}. Theorem~\ref{T:contour} says that on the event that the maximal trajectory has a zero at the origin, a Toom contour must be present. As stated precisely in Theorem~\ref{T:ToomPei}, this allows one to estimate $1-\ov\rho(p)$ from above by the expected number of Toom contours that are present in a cellular automaton. This is the core of Toom's Peierls argument. Theorem~\ref{T:cycle} shows that for a subclass of monotone cellular automata, it is possible to work with \emph{Toom cycles} which are Toom contours with additional useful properties that often lead to sharper bounds.

The remainder of the paper is devoted to proofs. In Section~\ref{S:prelim} we prove some preparatory results, in Section~\ref{S:Toomconstr} we prove the results from Section~\ref{S:ToPe} about Toom contours, and in Section~\ref{S:Peierls} we apply these results to prove stability of a class of monotone cellular automata with intrinsic randomness and derive some explicit bounds on the critical noise parameter.

\subsection{Stability of monotone cellular automata}\label{S:stab}

In this subsection we state a theorem giving sufficient conditions for the stability of monotone cellular automata with intrinsic randomness. The statement of the theorem involves edge speeds and gives additional insight into Toom's stability theorem, which it generalises.

\begin{defi}[Edge speed]\label{def:speed}
Let $\ell:\R^d\to\R$ be a linear form that is not identically zero and let $\phi:\{0,1\}^{\Z^d}\to\{0,1\}$ be a monotone local map. We call the quantity
\be\label{edge}
\eps_\phi(\ell):=\sup_{A\in\Oi(\phi)}\inf_{i\in A}\ell(i).
\ee
the \emph{edge speed} of $\phi$ in the direction $\ell$.
\end{defi}

The name ``edge speed'' already suggests its interpretation. For any linear form $\ell:\R^d\to\R$, let $H^\ell_r\in\{0,1\}^{\Z^d}$ denote the half-space configuration defined by
\be\label{halfspace}
H^\ell_r(i):=\left\{\ba{ll}
1\quad&\mbox{if }\ell(i)\geq r,\\[5pt]
0\quad&\mbox{if }\ell(i)<r
\ea\right. \quad (r\in \R).
\ee
The following lemma explains the name ``edge speed''.

\bl[Edge speeds]
Let\label{L:edge} $\ell:\R^d\to\R$ be a linear form that is not identically equal to zero and $\phi:\{0,1\}^{\Z^d}\to\{0,1\}$ be a monotone local map. Then for each $r\in\R$ and $t\geq 0$ the map from (\ref{Psiphi}) satisfies
\be\label{PsiH}
\Psi^t_\phi(H^\ell_r)=H^\ell_{r-t\eps_\phi(\ell)}.
\ee
\el

The result above, which follows easily from the definitions, is for completeness proved in Subsection~\ref{S:erod}. To state our stability result, we need one more definition.

\begin{defi}[Polar functions] Given\label{defi:polar} an integer $\sigma\geq 2$,
a \emph{polar function} of dimension $\sigma$ is a linear function
\be
\R^d\ni z\mapsto L(z)=(L_1(z),\ldots,L_\sig(z))\in\R^\sig
\ee
such that
\be\label{polar}
\sum_{s=1}^\sig L_s(z)=0\qquad(z\in\R^d).
\ee
\end{defi}

Note that a polar function can be regarded as a collection of $\sig$ points in $\mathbb{R}^d$ whose sum is zero.

In Subsection~\ref{S:mainproof} we will use Toom contours to prove the following stability result.

\bt[Stability of monotone cellular automata with intrinsic randomness]
Fix\label{T:detspeed} $m\geq 1$ and let $\Phi^0$ be a monotone cellular automaton defined by maps $\phi_1,\ldots,\phi_m$ and a probability distribution $\rbf(1),\ldots,\rbf(m)$.  Assume that there exists a linear polar function $L$ of dimension $\sig\geq 2$ such that the worst-case edge speeds
\be\label{detspeed}
\eps_s:=\inf_{1\leq k\leq m}\eps_{\phi_k}(L_s)\qquad(\lis)
\ee
satisfy 
\be\label{eps}
\eps:=\sum_{s=1}^\sig\eps_s>0.
\ee
Then $\Phi^0$ is stable.
\et

Theorem~\ref{T:detspeed} is far from optimal in terms of what can be achieved by Toom's Peierls argument, but to improve on it one needs more advanced methods to estimate the Peierls sum which will be presented in our companion paper \cite{SSTcont}. Although Theorem~\ref{T:detspeed} is suboptimal in the presence of intrinsic randomness, it is optimal in the deterministic case. To see this, we need the following alternative erosion criterion originally due to \cite[Lemma~12]{Pon13}, the proof of which will be given in Subsection~\ref{S:erod}.

\bl[Alternative erosion criterion]
Let\label{L:erode} $\phi:\{0,1\}^{\Z^d}\to\{0,1\}$ be a non-constant monotone function. Then $\phi$ is an eroder if and only if there exists a polar function $L$ of dimension $\sig\geq 2$ such that the edge speeds defined in (\ref{edge}) satisfy
\be\label{erode}
\sum_{s=1}^\sig\eps_{\phi}(L_s)>0.
\ee
\el

It is instructive to see why (\ref{erode}) implies that $\phi$ is an eroder. Given a configuration $x\in\{0,1\}^{\Z^d}$ containing finitely many zeros, let  the \emph{extent} of $x$ be defined as 
\be\label{extent}
{\rm ext}(x):=
\begin{cases}
\dis\sum_{s=1}^\sig r_s(x)\quad\mbox{with}\quad
r_s(x):=\sup\big\{L_s(i):i\in\Z^d,\ x(i)=0\big\} & \mbox{ if }  x\neq \un 1,\\
\dis-\infty & \mbox{ if } x= \un 1
\end{cases}
\ee
By the defining property (\ref{polar}) of a linear polar function, ${\rm ext}(x)\geq 0$ for each $x\neq\un 1$. Lemma~\ref{L:edge} and the monotonicity of $\phi$ imply that for each configuration $x$ with finitely many zeros,
\be
{\rm ext}(\Psi^t_\phi(x))\leq{\rm ext}(x)-\eps_\phi(L)t\qquad(t\geq 0),
\ee
and hence $\Psi^t_\phi(x)=\un 1$ for all $t\geq{\rm ext}(x)/\eps$. In the case with intrinsic randomness, condition (\ref{eps}) similarly implies that if $(X_t)_{t\geq 0}$ is the Markov chain defined as in (\ref{Markov}) in terms of the unperturbed automaton $\Phi^0$, started in an initial state $X_0=x$ with finitely many zeros, then almost surely
\be
{\rm ext}(X_t)\leq{\rm ext}(x)-\eps t\qquad(t\geq 0),
\ee
and $X_t=\un 1$ for all $t\geq{\rm ext}(x)/\eps$. Thus Theorem~\ref{T:detspeed} proves stability under the assumption that under the unperturbed evolution, finite collections of zeros disappear after a finite \emph{deterministic} time. There are many examples of monotone cellular automata with intrinsic randomness that do not satisfy (\ref{eps}) but for which under the unperturbed evolution, finite collections of zeros disappear after a finite \emph{random} time. For some of these, we will prove stability in our companion paper \cite{SSTcont}.

\subsection{Bounds on the critical noise parameter}\label{S:pc}

In this subsection we apply Theorem~\ref{T:detspeed} to some concrete examples and derive explicit bounds on the critical noise parameter $p_{\rm c}$ from (\ref{pc}).

We first set $m=1$ and consider the deterministic cellular automaton on $\Z^2$ defined by the single map $\phi_1=\phi^{\rm coop}$. The function $L:\R^2\to\R^2$ defined as
\be\label{eq:cooppolar}
L_1(z):=-z_1-z_2, \quad L_2(z):=z_1+z_2\qquad\big(z=(z_1,z_2)\in\R^2\big)
\ee
is a linear polar function of dimension $\sigma=2$ in the sense of Definition~\ref{defi:polar}. Using \eqref{ANEC} and \eqref{edge}, we see that the corresponding edge speeds from \eqref{detspeed} are given by
\be
\eps_1=\eps_{\phi^{\rm coop}}(L_1)=0, \quad \eps_2=\eps_{\phi^{\rm coop}}(L_2)=1,
\ee
so $\eps=\eps_1+\eps_2>0$ and hence Theorem~\ref{T:detspeed} implies that this cellular automaton is stable.

We next set $m=1$ and $\phi_1=\phi^{\rm NEC}$. We define a linear polar function $L$ of dimension $\sig=3$ by
\be\label{eq:Toompolar}
L_1(z_1,z_2):=-z_1,\quad
L_2(z_1,z_2):=-z_2,\quad
L_3(z_1,z_2):=z_1+z_2\qquad(z\in\R^2).
\ee
One can check that for this choice of $L$ (recall \eqref{ANEC})
\be
\eps_1=\eps_{\phi^{\rm NEC}}(L_1)=0, \quad \eps_2=\eps_{\phi^{\rm NEC}}(L_2)=0,
\quad \eps_3=\eps_{\phi^{\rm NEC}}(L_3)=1,
\ee
which implies $\eps=\eps_1+\eps_2+\eps_3>0$, hence stability.

To also see an example with intrinsic randomness, consider the case $m=2$ with $\phi_1=\phi^{\rm NEC}$ and $\phi_2=\phi^{\rm coop}$. Using the polar function \eqref{eq:Toompolar} one can check that
\begin{align}
\eps_1=\eps_{\phi^{\rm NEC}}(L_1)\wedge\eps_{\phi^{\rm coop}}(L_1)=0\wedge 0=0,\\
\eps_2=\eps_{\phi^{\rm NEC}}(L_2)\wedge\eps_{\phi^{\rm coop}}(L_2)=0\wedge 0=0,\\
\eps_3=\eps_{\phi^{\rm NEC}}(L_3)\wedge\eps_{\phi^{\rm coop}}(L_3)=1\wedge 1=1,
\end{align}
which implies $\eps>0$. Therefore, Theorem~\ref{T:detspeed} implies stability for this cellular automaton regardless of the choice of the probability distribution $\rbf=\big(\rbf(1),\rbf(2)\big)$ on $\{1,2\}$.

To see an example where Theorem~\ref{T:detspeed} is not applicable, consider the case $m=4$ with $\phi_1=\phi^{\rm NEC}$, $\phi_2=\phi^{\rm NWC}$, $\phi_3=\phi^{\rm SWC}$, and $\phi_4=\phi^{\rm SEC}$. In this case, there exists no polar function that satisfies the hypothesis of Theorem~\ref{T:detspeed}. If $\rbf$ is the uniform distribution on $\{1,2,3,4\}$, then there is no direction in which the ones tend to invade the zeros or vice versa. In other words, the ``effective edge speed'' is zero in each direction. On the other hand, if $\rbf(1)\neq\rbf(3)$ or $\rbf(2)\neq\rbf(4)$, then it seems likely there is a positive ``effective edge speed''. We make the following precise conjecture.

\bcon[Random direction NEC voting]
The\label{C:symNEC} cellular automaton that applies the maps $\phi^{\rm NEC}$, $\phi^{\rm NWC}$, $\phi^{\rm SWC}$, and $\phi^{\rm SEC}$ with probabilities $\rbf(1)$, $\rbf(2)$, $\rbf(3)$, and $\rbf(4)$, respectively, is stable if $\rbf(1)\neq\rbf(3)$ or $\rbf(2)\neq\rbf(4)$, and unstable if $\rbf(1)=\rbf(3)$ and $\rbf(2)=\rbf(4)$.
\econ

The proof of Theorem~\ref{T:detspeed} allows us to derive explicit lower bounds on the critical noise parameter $p_{\rm c}$ from (\ref{pc}). In particular, in Subsection~\ref{S:explic} we will prove the following bounds.

\bp[Explicit bounds] 
For\label{P:coopbd} the deterministic cellular automaton on $\Z^2$ that applies $\phi^{\rm coop}$ in each space-time point $p_{\rm c}\geq 1/64$. For the deterministic cellular automaton on $\Z^2$ that applies $\phi^{\rm NEC}$ in each space-time point $p_{\rm c}\geq 3^{-21}$.
\ep

In our companion paper \cite{SSTcont}, using a more advanced method to bound the Peierls sum, we will improve the lower bound for the cellular automaton defined by $\phi^{\rm NEC}$ to $p_{\rm c}\geq 1/12000$. Numerical simulations suggest that the true value of  $p_{\rm c}$ is $\approx 0.105$ for $\phi^{\rm coop}$ and $\approx 0.053$  for $\phi^{\rm NEC}$. There is a good reason why the rigorous bounds for $\phi^{\rm NEC}$ are worse than for $\phi^{\rm coop}$. If we want to apply Lemma~\ref{L:erode} to prove that $\phi^{\rm NEC}$ is an eroder, then we need a linear polar function of dimension at least three, while for $\phi^{\rm coop}$ a linear polar function of dimension two suffices. In general, the higher the dimension of the linear polar function, the worse the bounds. For linear polar functions of dimension two, we can moreover use Toom cycles instead of Toom contours, which also leads to sharper bounds.

\subsection{History of the problem and motivation}\label{S:discussion}

The cellular automaton defined by the NEC voting map $\phi^{\rm NEC}$ is nowadays known as \emph{Toom's model}. In line with Stigler's law of eponymy, Toom's model was not invented by Toom, but by Vasilyev, Petrovskaya, and Pyatetski-Shapiro, who simulated random perturbations of this and other models on a computer \cite{VPP69}. Toom, having heard of \cite{VPP69} during a seminar, proved in \cite{Too74} that there exist random cellular automata on $\Z^d$ with at least $d$ different invariant laws. Although Toom's model is not explicitly mentioned in the paper, his proof method can be applied to prove that $p_{\rm c}>0$ for his model. In \cite{Too80}, Toom improved his methods and proved his celebrated stability theorem. His paper is quite hard to read. A more accessible account of Toom's original argument (with pictures!) in the special case of Toom's model can be found in the appendix of \cite{LMS90}.\footnote{Unfortunately, their Figure~6 contains a small mistake, in the form of an arrow that should not be there.}

Bramson and Gray \cite{BG91} have given another alternative proof of Toom's stability theorem that relies on comparison with continuum models (which describe unions of convex sets in $\R^d$ evolving in continuous time) and renormalisation-style block arguments. A disadvantage of this approach is that it is restricted to lattices that can be rescaled to $\R^d$ while Toom's method can also work on lattices such as trees, as demonstrated in \cite{CSS24}. Gray \cite{Gra99} proved a stability theorem for monotone interacting particle systems (i.e., in continuous time). The proofs use ideas from \cite{Too80} and \cite{BG91}. Gray also derived necessary and sufficient conditions for a monotonic map to be an eroder \cite[Thm~18.2.1]{Gra99}, apparently overlooking the fact that Toom had already proved the much simpler condition (\ref{erosion}).

The cellular automaton that applies the monotone map $\phi$ with probability $p$ and the identity map $\phi^{\rm id}$ with probability $1-p$ is also referred to in the literature as $p$-\textit{asynchronous} cellular automaton. In asynchronous cellular automata, cells do not update their states simultaneously. There are various ways to define this asynchrony; for a comprehensive survey, see \cite{Fat13}. In \cite{Gha92}, a generalisation of Toom's theorem was presented for a particular class of asynchronous cellular automata.

Motivated by abstract problems in computer science, a number of authors have given alternative proofs of Toom's stability theorem in a more restrictive setting \cite{GR88,BS88,Gac95,Gac21}. Their main interest is in a three-dimensional system which evolves in two steps: letting $e_1,e_2,e_3$ denote the basis vectors in $\Z^3$, they first replace $X_n(i)$ by
\[
X'_n(i):={\tt round}\big((X_n(i)+X_n(i+e_1)+X_n(i+e_2))/3\big),
\]
and then set
\[
X_{n+1}(i):={\tt round}\big((X'_n(i)+X'_n(i+e_3)+X'_n(i-e_3))/3\big).
\]
They prove explicit bounds for finite systems, although for values of $p$ that are extremely close to zero.\footnote{In particular, \cite{Gac95} needs $p<2^{-21}3^{-8}$.} The proofs of \cite{GR88} do not use Toom's Peierls argument but rely on different methods. Their bounds were improved in \cite{BS88}. Still better bounds can be found in the unpublished note \cite{Gac95}. The proofs in the latter manuscript are very similar to Toom's argument, with some crucial suggested improvements at the end that are hard to follow due to missing definitions. This version of the argument seems to have inspired the incomplete note by John Preskill \cite{Pre07} who links it to the interesting idea of counting ``minimal explanations''. We will use this general idea in Subsection~\ref{S:minexpl} below, but our precise definition of a ``minimal explanation'' differs a bit from his. As explained at Figure~\ref{fig:minexpl1} and at the end of Subsection~\ref{S:presence}, the relation between Toom contours and minimal explanations is not so straightforward as suggested in \cite{Gac95,Pre07}.

Around 1985, Durrett and Gray submitted a very interesting paper about an interacting particle system based on the map $\phi^{\rm coop}$ from (\ref{phiNEC}). The major revision requested by the referee never materialised, however. For many years, a short note announcing the results without proofs \cite{Dur86} was the only accessible source to this material but recently Rick Durrett has made the original preprint available on his homepage \cite{DG85}. Hwa-Nien Chen \cite{Che92,Che94}, who was a PhD student of Lawrence Gray, studied the stability of various variations of Toom's model under perturbations of the initial state and the birth rate. The proofs of two of his four theorems depend on results that he cites from the preprint \cite{DG85}. Ponselet \cite{Pon13} gave an excellent account of the existing literature and together with her supervisor proved exponential decay of correlations for the upper invariant law of a large class of randomly perturbed monotone cellular automata \cite{MP11}.

There exists duality theory for general monotone interacting particle systems \cite{Gra86,SS18,LS23}. The basic idea is that the state in the origin at time zero is a monotone function of the state at time $-t$, and this monotone function evolves in a Markovian way as a function of~$t$. As noted in \cite{Dur86} this dual process plays an important ingredient of the proofs of \cite{DG85}. It is also closely related to the minimal explanations of Preskill \cite{Pre07}. A good understanding of this dual process could potentially help solve many open problems in the area, but its behaviour is already quite complicated in the mean-field case \cite{MSS20}.

\begin{figure}[htb]
\begin{center}
\includegraphics[width=\textwidth]{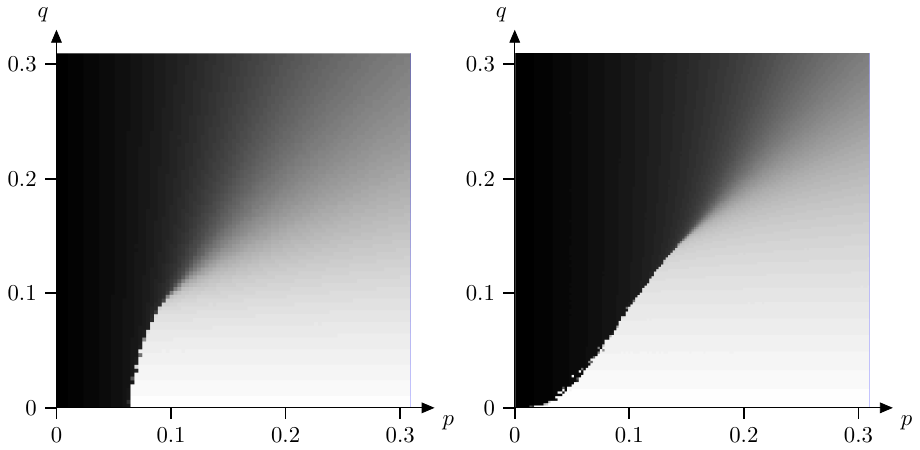}
\caption{Density $\ov\rho$ of the upper invariant law of two monotone random cellular automata as a function of the parameters, shown on a scale from 0 (white) to 1 (black). On the left: a version of Toom's model that applies the maps $\phi^0$, $\phi^1$, and $\phi^{\rm NEC}$ with probabilities $p$, $q$, and $1-p-q$, respectively. On the right: the mononotone random cellular automaton that applies the maps $\phi^0$, $\phi^1$, and $\phi^{\rm NN}$ with probabilities $p$, $q$, and $1-p-q$, respectively. The map $\phi^{\rm NEC}$ is an eroder but $\phi^{\rm NN}$ is not. By the symmetry between the 0's and the 1's, in both models, the density $\un\rho(p, q)$ of the lower invariant law equals $1-\ov\rho(q,p)$. Due to metastability effects, the area where the upper invariant law differs from the lower invariant law is shown too large in these numerical data. For Toom's model with $q=0$, the data shown above suggest a first order phase transition at $p_{\rm c}\approx 0.057$ but based on numerical data for edge speeds we believe the true value is $p_{\rm c}\approx 0.053$. We conjecture that the model on the right has a unique invariant law everywhere except on the diagonal $p=q$ for $p$ sufficiently small. (The pictures are based on a combination of simulations of various precision. In the most sensitive parameter regime, the system size is $100\times 100$ with periodic boundary conditions and the total time is $1000$.)}
\label{fig:updens}
\end{center}
\end{figure}

In numerical investigations of monotone cellular automata, it is often useful to take a wider view and perturb the system not only with i.i.d.\ zeros but also with i.i.d.\ ones. Recall the constant maps $\phi^0$ and $\phi^1$ defined in (\ref{consmap}). If $\phi:\{0,1\}^{\Z^d}\to\{0,1\}$ is a monotone local map that is not constant, then it is interesting to look at i.i.d.\ collections of random variables $\Phi^{p,q}=(\Phi^{p,q}_{i,t})_{(i,t)\in\Z^{d+1}}$ with values in $\{\phi^0,\phi^1,\phi\}$ such that
\be
\P\big[\Phi^{p,q}_{i,t}=\phi^0\big]=p,\quad\P\big[\Phi^{p,q}_{i,t}=\phi^1\big]=q,\quad\P\big[\Phi^{p,q}_{i,t}=\phi\big]=1-p-q,
\ee
where $p,q\geq 0$ with $p+q\leq 1$ are parameters of the model. Figure~\ref{fig:updens} shows numerical data for the density $\ov\rho(p,q)$ of the upper invariant law for such a cellular automaton in the case when $\phi=\phi^{\rm NEC}$ and $\phi=\phi^{\rm NN}$, respectively. We see that in line with Toom's theorem, $\lim_{p\to 0}\ov\rho(p,0)=1$ for the NEC voting rule but $\ov\rho(p,0)=0$ for all $p>0$ in the case of NN voting. Nevertheless, the simulations suggest that the NN voting rule is stable under symmetric noise, in the sense that $\lim_{p\to 0}\ov\rho(p,p)=1$. Proving this is a long-standing open problem; a continuous-time version of this model is mentioned in \cite[Example~I.4.3(e)]{Lig85}. A closely related problem, that is also open, is to show that for the NEC voting rule the function $p\mapsto\rho(p)$ makes a jump at $p_{\rm c}:=\sup\{p:\ov\rho(p)>0\}$.

\section{Toom contours}\label{S:ToPe}

\subsection{Monotone cellular automata}\label{S:monaut}

In this section, we introduce Toom contours, which are the central object in Toom's Peierls argument. Toom contours can be defined for monotone cellular automata in which space-time has a more general structure than $\Z^{d+1}$. In the present subsection, we extend the definitions of Section~\ref{S:intmain} to this more general set-up.

Let $\La$ be a countable set. As in Subsection~\ref{S:intro}, we say that a map $\phh:\{0,1\}^\La\to\{0,1\}$ is \emph{local} if there exists a finite $\De\sub\La$ such that $\phh(x)$ depends only on $(x(i))_{i\in\De}$ and we let $\De(\phh)$ denote the smallest such set. In analogy with (\ref{consmap}), in the present setting, we denote the constant functions by 
\be\label{phh01}
\phh^0(x):=0\quand\phh^1(x):=1\qquad(x\in\{0,1\}^\La).
\ee
We let $\Oi(\phh)$ denote the set of minimal one-sets of $\phh$, defined as in Subsection \ref{S:intro}.

Recall that a \emph{directed graph} is a pair $(V,\vec E)$ where $V$ is a set whose elements are called \emph{vertices} and $\vec E$ is a subset of $V\times V$ whose elements are called \emph{directed edges}. For each directed edge $(v,w)\in\vec E$, we call $v$ the starting vertex and $w$ the endvertex. We say that $(V,\vec E)$ is \emph{acyclic} if there do not exist $n\geq 1$ and $v_0,\ldots,v_n\in V$ with $v_n=v_0$ such that $(v_{k-1},v_k)\in\vec E$ for all $0<k\leq n$.

Let $\phhb=(\phh_i)_{i\in\La}$ be a collection of local maps $\phh_i:\{0,1\}^\La\to\{0,1\}$, and let
\be\label{vecF}
\vec H(\phhb):=\big\{(i,j)\in\La^2:j\in\De(\phh_i)\big\}.
\ee
Then $(\La,\vec H(\phhb))$ is a directed graph. Generalising our earlier definition, we define a \emph{cellular automaton} to be a collection of local maps $\phhb=(\phh_i)_{i\in\La}$ for which the directed graph $(\La,\vec H(\phhb))$ is acyclic. We call $(\La,\vec H(\phhb))$ the \emph{dependence graph} associated with $\phhb$. A \emph{trajectory} of a cellular automaton is a function $x:\La\to\{0,1\}$ such that
\be\label{traj}
x(i)=\phh_i(x)\qquad(i\in\La).
\ee
We can think of $x$ as a function that describes the state of a cellular automaton as a function of space-time. Then (\ref{traj}) says that the state in the space-time point $i$ is a function of the states in the space-time points from the set $\De(\phi_i)$, which we think of as preceding $i$. A cellular automaton is \emph{monotone} if $\phh_i$ is a monotone map for each $i\in\La$, i.e., $x\leq y$ (coordinatewise) implies $\phh_i(x)\leq\phh_i(y)$.

To make the link with our earlier definitions from Subsection~\ref{S:intro}, let $\Phi^p=(\Phi^p_{i,t})_{(i,t)\in\Z^{d+1}}$ be a monotone cellular automaton of the type considered in (\ref{Phip}), and for each $(i,t)\in\Z^{d+1}$, define $\Phi^p_{(i,t)}:\{0,1\}^{\Z^{d+1}}\to\{0,1\}$ by
\be\label{Phip2}
\Phi^p_{(i,t)}\big((x(i',t'))_{(i',t')\in\Z^{d+1}}\big)
:=\Phi^p_{i,t}\big((x(i+i',t-1))_{i'\in\Z^d}\big)
\qquad\big(x\in\{0,1\}^{\Z^{d+1}}\big).
\ee
Then $(\Phi^p_{(i,t)})_{(i,t)\in\Z^{d+1}}$ is a random monotone cellular automaton according to the definitions of the present section. Note the subtle difference in notation between $\Phi^p_{i,t}$ and $\Phi^p_{(i,t)}$. By a slight abuse of notation, we use the symbol $\Phi^p$ for both the collections $(\Phi^p_{i,t})_{(i,t)\in\Z^{d+1}}$ and $(\Phi^p_{(i,t)})_{(i,t)\in\Z^{d+1}}$.

We next turn our attention to the lower and upper invariant laws from formula (\ref{toupper}). The following two lemmas introduce two closely related objects, the minimal and maximal trajectories, and show how they are related to the lower and upper invariant laws. We prove these lemmas in Subsection~\ref{S:maxtraj}.

\bl[Minimal and maximal trajectories]
Let\label{L:maxtraj} $\phhb$ be a monotone cellular automaton. Then there exist trajectories $\un x$ and $\ov x$ that are uniquely characterised by the property that each trajectory $x$ of $\phhb$ satisfies $\un x\leq x\leq\ov x$ (pointwise).
\el

\bl[Lower and upper invariant laws]
Let\label{L:numax} $\Phi^p$ be the random monotone cellular automaton defined in (\ref{Phip2}) and let $\un X^p$ and $\ov X^p$ be its minimal and maximal trajectories. Then
\be
\P\big[\big(\un X^p(i,t)\big)_{i\in\Z^d}\in\,\cdot\,\big]=\un\nu_p
\quand
\P\big[\big(\ov X^p(i,t)\big)_{i\in\Z^d}\in\,\cdot\,\big]=\ov\nu_p
\qquad(t\in\Z),
\ee
where $\un\nu_p$ and $\ov\nu_p$ are the lower and upper invariant laws of the Markov chain in (\ref{Markov}).
\el

Let $\Phi$ be a random monotone cellular automaton, i.e., a random variable taking values in the space of all monotone cellular automata on a given space-time set $\La$, and let $\ov X$ denote its maximal trajectory, which is now also random. In Theorem~\ref{T:ToomPei} below, we give a lower bound on the probability $\P[\ov X(i)=1]$. We will show that on the event that $\ov X(i)=0$, the random monotone cellular automaton $\Phi$ must contain a certain structure that we will call a \emph{Toom contour rooted at} $i$. The probability that $\ov X(i)=0$ can then be estimated from above by the expected number of Toom contours rooted at $i$ that are present in $\Phi$. In particular, applying this to the maximal trajectory $\ov X^p$ of the random monotone cellular automaton $\Phi^p$, we are under certain additional assumptions able to show that the density $\ov\rho(p)=\P[\ov X^p(i,t)=1]$ of the upper invariant law, which in this case does not depend on $(i,t)$, tends to one as $p\to 0$. In its essence, the method goes back to Toom's proof of \cite[Thm~5]{Too80} but we have significantly modified and simplified the argument with the aim of making it more flexible and intuitive. At the end of Subsection~\ref{S:presence} we give an overview of the most significant differences between our formulation of Toom's Peierls argument and the original formulation in \cite{Too80}.

\subsection{Toom contours}\label{S:con}

We will need directed graphs in which both the vertices and the edges can have different types. Let $A$ and $B$ be finite sets. By definition, a \emph{typed directed graph} with \emph{vertex set} $V$, \emph{vertex type set} $A$, and \emph{edge type set} $B$ is a pair $(\Vi,\Ei)$ where $\Vi$ is a subset of $V\times A$ and $\Ei$ is a subset of $V\times V\times B$, such that
\be\label{typedef}
\forall v\in V\ \exists a\in A\mbox{ s.t.\ }(v,a)\in\Vi.
\ee
For each $a\in A$ and $b\in B$, we call
\be
V_a:=\big\{v:(v,a)\in\Vi\big\}\quand
\vec E_b:=\big\{(v,w):(v,w,b)\in\Ei\big\}
\ee
the set of vertices of type $a$ and the set of directed edges of type $b$, respectively. Note that vertices can have more than one type, i.e., $V_a$ and $V_{a'}$ are not necessarily disjoint for $a\neq a'$, and the same applies to edges. As a consequence, several edges of different types can connect the same two vertices $v,w$, but always at most one of each type. If $(\Vi,\Ei)$ is a typed directed graph, then we let $(V,\vec E)$ denote the directed graph given by
\be
V=\bigcup_{a\in A}V_a\quand\vec E:=\bigcup_{b\in B}\vec E_b,
\ee
where the first equality follows from (\ref{typedef}) and the second equality is a definition. We call $(V,\vec E)$ the \emph{untyped} directed graph associated with $(\Vi,\Ei)$. We also set $E:=\big\{\{v,w\}:(v,w)\in\vec E\big\}$. Then $(V,E)$ is an undirected graph, which we call the undirected graph \emph{associated with} $(V,\vec E)$. We say that a typed directed graph $(\Vi,\Ei)$ or a directed graph $(V,\vec E)$ are \emph{connected} if their associated undirected graph $(V,E)$ is connected. A \emph{rooted} directed graph is a triple $(v_\circ,V,\vec E)$ such that $(V,\vec E)$ is a directed graph and $v_\circ\in V$ is a specially designated vertex, called the \emph{root}. Rooted undirected graphs and rooted typed directed graphs are defined in the same way.

For any directed graph $(V,\vec E)$, we let
\be
\vec E_{\rm in}(v):=\big\{(u,v')\in\vec E:v'=v\big\}
\quand
\vec E_{\rm out}(v):=\big\{(v',w)\in\vec E:v'=v\big\}
\ee
denote the sets of directed edges entering and leaving a given vertex $v\in V$, respectively. Similarly, in a typed directed graph, $\vec E_{b,{\rm in}}(v)$ and $\vec E_{b,{\rm out}}(v)$ denote the sets of incoming or outgoing directed edges of type~$b$ at $v$.

We adopt the following general notation. For any directed graph $(V,\vec E)$, set $\La$, and function $\psi:V\to\La$, we let
\be\label{psiedge}
\psi(V):=\big\{\psi(v):v\in V\big\}\quand
\psi(\vec E):=\big\{\big(\psi(v),\psi(w)\big):(v,w)\in\vec E\big\}
\ee
denote the images of $V$ and $\vec E$ under $\psi$. We can naturally view $\big(\psi(V),\psi(\vec E)\big)$ as a directed graph with set of vertices $\psi(V)$ and set of directed edges $\psi(\vec E)$. We denote this graph by $\psi(V,\vec E):=\big(\psi(V),\psi(\vec E)\big)$. Similarly, if $(\Vi,\Ei)$ is a typed directed graph, then we let $\psi(\Vi,\Ei)$ denote the typed directed graph defined as
\be\ba{r}\label{psiEi}
\dis\psi(\Vi,\Ei):=\big(\psi(\Vi),\psi(\Ei)\big)
\qquad\quad\mbox{with}\quad
\psi(\Vi):=\big\{\big(\psi(v),a\big):(v,a)\in\Vi\big\}\\[5pt]
\dis\quand
\psi(\Ei):=\big\{\big(\psi(v),\psi(w),b\big):(v,w,b)\in\Ei\big\}.
\ec
Also, if $(v_\circ,V,\vec E)$ is a rooted directed graph, then we let $\psi(v_\circ,V,\vec E)$ denote the rooted directed graph $\big(\psi(v_\circ),\psi(V),\psi(\vec E)\big)$, and we use similar notation for rooted typed directed graphs. Two typed directed graphs $(\Vi,\Ei)$ and $(\Wi,\Fi)$ are \emph{isomorphic} if there exists a bijection $\psi:V\to W$ such that $\psi(\Vi,\Ei)=(\Wi,\Fi)$. Similar conventions apply to directed graphs, rooted directed graphs, and so on.

\begin{figure}[t]
\begin{center}
\includegraphics[width=0.7\textwidth]{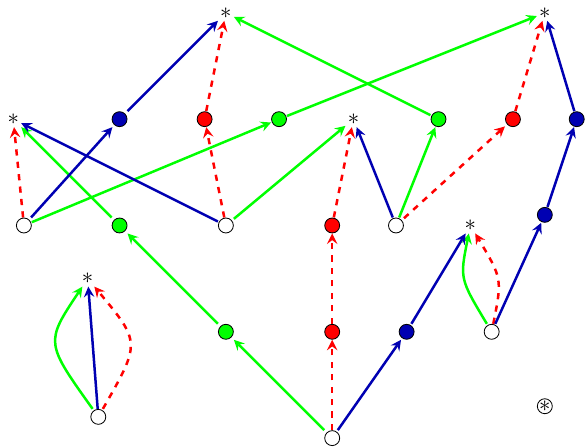}
\caption{Example of a Toom graph with three charges. Sources are indicated with open dots, sinks with asterisks, and internal vertices and edges of the three possible charges with three colours. Note the isolated vertex in the lower right corner, which is a source and a sink at the same time.}
\label{fig:Toomgraph}
\end{center}
\end{figure}

\begin{defi}\label{D:toomgraph}
A \emph{Toom graph} with $\sig\geq 1$ \emph{charges} is a typed directed graph $(\Vi,\Ei)$ with vertex type set $\{\circ,\ast,1,\ldots,\sig\}$ and edge type set $\{1,\ldots,\sig\}$ that satisfies the following conditions:
\begin{enumerate}
\item $|\vec E_{s,{\rm in}}(v)|=0$ $(\lis)$ and $|\vec E_{1,{\rm out}}(v)|=\cdots=|\vec E_{\sig,{\rm out}}(v)|\leq 1$ for all $v\in V_\circ$.
\item $|\vec E_{s,{\rm out}}(v)|=0$ $(\lis)$ and $|\vec E_{1,{\rm in}}(v)|=\cdots=|\vec E_{\sig,{\rm in}}(v)|\leq 1$ for all $v\in V_\ast$.
\item $|\vec E_{s,{\rm in}}(v)|=1=|\vec E_{s,{\rm out}}(v)|$ and $|\vec E_{l,{\rm in}}(v)|=0=|\vec E_{l,{\rm out}}(v)|$ for each $l\neq s$ and $v\in V_s$.
\end{enumerate}
\end{defi}

See Figure~\ref{fig:Toomgraph} for a picture of a Toom graph with three charges. Vertices in $V_\circ,V_\ast$, and $V_s$ are called \emph{sources}, \emph{sinks}, and \emph{internal vertices} with \emph{charge} $s$, respectively. Vertices in $V_\circ\cap V_\ast$ are called \emph{isolated vertices}. With the exception of isolated vertices, the inequalities $\leq 1$ in (i) and (ii) are equalities. Informally, we can imagine that at each source there emerge $\sig$ charges, one of each type, that then travel via internal vertices of the corresponding charge through the graph until they arrive at a sink, in such a way that at each sink there converge precisely $\sig$ charges, one of each type. This informal picture holds even for isolated vertices, if we imagine that in this case, the charges arrive immediately at the sink that is at the same time a source. It is clear from this informal picture that $|V_\circ|=|V_\ast|$, i.e., the number of sources equals the number of sinks. We let $(V,\vec E)$ denote the directed graph associated with $(\Vi,\Ei)$.

Toom graphs and the Toom contours that will be defined below were designed to make use of the condition (\ref{eps}) of Theorem~\ref{T:detspeed} on the worst-case edge speeds. The curious reader may skip ahead to the beginning of Subsection~\ref{S:edgebnd} where we give an informal description of the main idea of the proof of Theorem~\ref{T:detspeed}.

Recall that a rooted directed graph is a directed graph with a specially designated vertex, called the root. In the case of Toom graphs, we will always assume that the root is a source.

\begin{defi}\label{D:toomroot}
A \emph{rooted Toom graph} with $\sig\geq 1$ \emph{charges} is a rooted typed directed graph $(v_\circ,\Vi,\Ei)$ such that $(\Vi,\Ei)$ is a Toom graph with $\sig\geq 1$ charges and $v_\circ\in V_\circ$. For any rooted Toom graph $(v_\circ,\Vi,\Ei)$, we write
\be\label{toomroot}
V'_\circ:=V_\circ\beh\{v_\circ\}\quand V'_s:=V_s\cup\{v_\circ\}\quad(\lis).
\ee
\end{defi}

The idea behind (\ref{toomroot}) is that for rooted Toom contours, we view the root more as if it were a collection of internal vertices than as a source. This is reflected in condition~(ii) of the following definition.

\begin{defi}\label{D:embed}
Let $(v_\circ,\Vi,\Ei)$ be a rooted Toom graph and let $\La$ be a countable set. An \emph{embedding} of $(v_\circ,\Vi,\Ei)$ in $\La$ is a map $\psi:V\to\La$ such that:
\begin{enumerate}
\item $\psi(v_1)\neq\psi(v_2)$ for each $v_1\in V_\ast$ and $v_2\in V$ with $v_1\neq v_2$,
\item $\psi(v_1)\neq\psi(v_2)$ for each $v_1,v_2\in V'_s$ with $v_1\neq v_2$ $(\lis)$.
\end{enumerate}
\end{defi}

Condition~(i) says that sinks do not overlap with other vertices and condition~(ii) says that internal vertices do not overlap with other internal vertices of the same charge, where in line with (\ref{toomroot}) we view the root as a collection of internal vertices. We make the following observation.

\bl[No double incoming edges]
Let\label{L:nodoub} $\psi$ be an embedding of a rooted Toom graph $(v_\circ,\Vi,\Ei)$ with $\sig\geq 1$ edges in a set $\La$. Then
\be
\big|\{(v,w)\in\vec E_s:\psi(w)=j\}\big|\leq 1\qquad(j\in\La,\ \lis).
\ee 
\el

\bpro
Immediate from Definition~\ref{D:embed}, since each charged edge ends in an internal vertex of the same charge or in a sink.
\epro

\begin{defi}\label{D:contour}
Let $\La$ be a countable set. A \emph{Toom contour} in $\La$ with $\sig\geq 1$ charges is a quadruple $(v_\circ,\Vi,\Ei,\psi)$, where $(v_\circ,\Vi,\Ei)$ is a rooted connected Toom graph with $\sig$ charges and $\psi$ is an embedding of $(v_\circ,\Vi,\Ei)$ in $\La$. We say that the Toom contour is \emph{rooted} at $i_\circ:=\psi(v_\circ)$. Two Toom contours $(v_\circ,\Vi,\Ei,\psi)$ and $(v'_\circ,\Vi',\Ei',\psi')$ are \emph{isomorphic} if there exists a bijection $\chi:V\to V'$ such that $\chi(v_\circ,\Vi,\Ei)=(v'_\circ,\Vi',\Ei')$ and $\psi(v)=\psi'(\chi(v))$ $(v\in V)$. We say that $(v_\circ,\Vi,\Ei,\psi)$ and $(v'_\circ,\Vi',\Ei',\psi')$ are \emph{equivalent} if, using notation introduced in (\ref{psiEi}), one has $\psi(\Vi,\Ei)=\psi'(\Vi',\Ei')$.
\end{defi}

\begin{figure}[htb]
\begin{center}
\includegraphics[width=\textwidth]{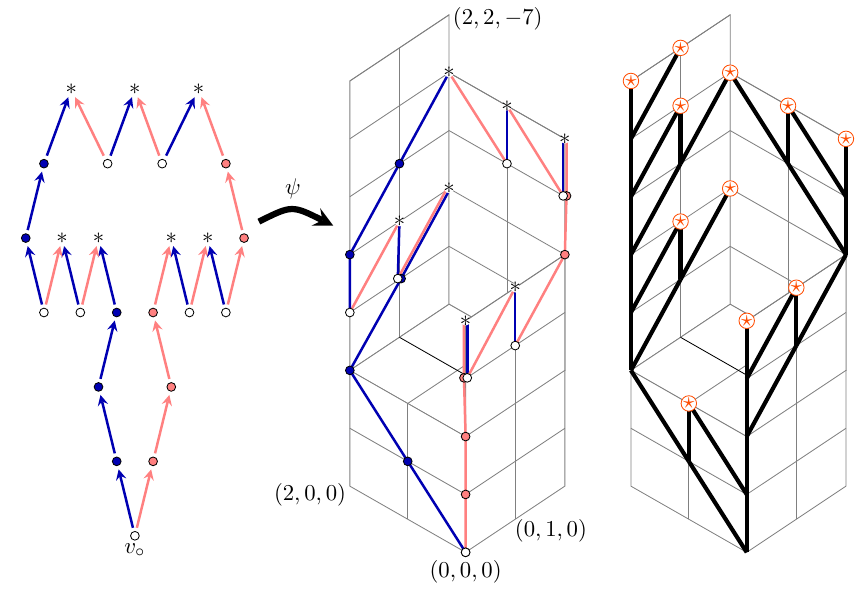}
\caption{A Toom contour in $\Z^3$ rooted at $(0,0,0)$. The third coordinate represents time and is plotted downwards. The picture on the right shows a minimal explanation (or rather its associated undirected explanation graph as defined in Subsection~\ref{S:minexpl}) for a monotone cellular automaton $\Phh^p$ that applies the maps $\phi^0$ and $\phi^{\rm coop}$ with probabilities $p$ and $1-p$, respectively. The origin has the value zero because the sites marked with a star are defective; removing any of these defective sites results in the origin having the value one. The Toom contour in the middle picture is present in $\Phh^p$. In particular, the sinks of the Toom contour coincide with some, though not with all of the defective sites of the minimal explanation.}
\label{fig:minexpl1}
\end{center}
\end{figure}

We note that as a result of Lemma~\ref{L:nodoub}, each charged edge in $\psi(\Ei)$ corresponds to a unique charged edge in $\Ei$. Two isomorphic Toom contours are clearly equivalent, but the converse implication does not hold, since sources can overlap with each other and with internal vertices and as a result, although two equivalent Toom contours have charged edges in the same locations, these edges can be differently connected leading to two Toom contours that are not isomorphic. See Figure~\ref{fig:minexpl1} for an example of a Toom contour with two charges.

\subsection{Presence of Toom contours}\label{S:presence}

Our next aim is to define when a Toom contour is present in a monotone cellular automaton $\phhb=(\phh_i)_{i\in\La}$. This will require us to make some extra assumptions and equip $\phhb$ with some extra structure.

\begin{defi}\label{D:typdep}
A \emph{typed dependence graph} with $\sig\geq 1$ types of edges is a typed directed graph $(\La,\Hi)$ with vertex type set $\{0,1,\bullet\}$ and edge type set $\{1,\ldots,\sig\}$ such that for $\vec H_s:=\{(i, j): (i, j, s)\in\Hi\}$
\begin{enumerate}
\item $\vec H_{s,{\rm out}}(i)=\emptyset$ for all $i\in\La_0\cup\La_1$ and $\lis$,
\item $\vec H_{s,{\rm out}}(i)\neq\emptyset$ for all $i\in\La_\bullet$ and $\lis$,
\end{enumerate}
and its associated untyped directed graph $(\La,\vec H)$ is acyclic. The monotone cellular automaton $\phhb=(\phh_i)_{i\in\La}$ \emph{associated with} the typed dependence graph $(\La,\Hi)$ is defined by
\be\label{phhFi}
\phh_i(x)=\left\{\ba{ll}
\dis\bigvee_{s=1}^\sig\bigwedge_{j:\,(i,j)\in\vec H_s}x(j)\quad&\dis\mbox{if }i\in\La_\bullet,\\[5pt]
\dis r\quad&\dis\mbox{if }i\in\La_r\quad(r=0,1),
\ea\right.
\ee
$\big(i\in\La,\ x\in\{0,1\}^\La\big)$.
\end{defi}

It is easy to see that if $(\La,\Hi)$ is a typed dependence graph, $\phhb$ is its associated monotone cellular automaton, and $(\La,\vec H)$ is its associated untyped directed graph, then $(\La,\vec H)$ is the dependence graph of $\phhb$ as defined in Subsection~\ref{S:monaut}. In particular, the assumption that $(\La,\vec H)$ is acyclic guarantees that (\ref{phhFi}) defines a cellular automaton. It is clear from (\ref{phhFi}) that $\phh_i$ is monotone for each $i\in\La$ and that $\phh_i$ is one of the constant maps $\phh^r$ $(r=0,1)$ defined in (\ref{phh01}) if and only if $i\in\La_r$ $(r=0,1)$. Furthermore, recalling~\eqref{Aphi} and the definition of one-sets, we can see that for each $i\in\La_\bullet$ the sets $\{j\in\Lambda: (i, j)\in\vec H_s\} \; (\lis)$ are one-sets of $\phh_i$, though not necessarily minimal. Elements of $\La_0$, where the constant zero map is applied, are called \emph{defective} sites. Below, we make use of the definition (\ref{toomroot}), i.e., we treat the root as if it were a collection of internal vertices.

\begin{defi}\label{D:present}
Let $(\La,\Hi)$ be a typed dependence graph with $\sig\geq 1$ types of edges. We say that a Toom contour $(v_\circ,\Vi,\Ei,\psi)$ with $\sig$ charges is \emph{present} in $(\La,\Hi)$ if:
\begin{enumerate}
\item $\dis\psi(v)\in\La_0$ for all $\dis v\in V_\ast$,
\item $\dis\big(\psi(v),\psi(w)\big)\in\vec H_s$ for all $(v,w)\in\vec E^\bullet_s$ $(\lis$),
\item $\dis\big(\psi(v),\psi(w)\big)\in\vec H$ for all $(v,w)\in\vec E^\circ$,
\end{enumerate}
where for any rooted Toom graph $(v_\circ,\Vi,\Ei)$, we write
\be\ba{l}\label{Ecirc}
\vec E^\bullet:=\bigcup_{s=1}^\sig\vec E^\bullet_s
\quad\mbox{with}\quad
\vec E^\bullet_s:=\big\{(v,w)\in\vec E_s:v\in V'_s\big\}
\quad(\lis),\\[5pt]
\vec E^\circ:=\bigcup_{s=1}^\sig\vec E^\circ_s
\quad\mbox{with}\quad
\vec E^\circ_s:=\big\{(v,w)\in\vec E_s:v\in V'_\circ\big\}
\quad(\lis).
\ec
\end{defi}

Condition (i) says that sinks of the Toom contour correspond to defective sites of the typed dependence graph. Conditions (ii) and (iii) say that directed edges of the Toom graph $(\Vi,\Ei)$ are mapped to directed edges of the typed dependence graph $(\La,\Hi)$, where edges coming out of an internal vertex must be mapped to edges of the corresponding type, and we treat the root as if it were a collection of internal vertices. Let $W:=\psi(V)$ and $W_\ast:=\psi(V_\ast)$. We note that Definition~\ref{D:present} implies that
\be\label{noLa1}
W\cap\La_0=W_\ast\quand W\cap\La_1=\emptyset.
\ee
Indeed, the inclusion $W_\ast\sub W\cap\La_0$ is immediate from condition~(i) while conditions (ii) and (iii) imply $W\beh W_\ast\sub\La\beh(\La_0\cup\La_1)$ since for $v\in V\beh V_\ast$, one has $\vec E_{\rm out}(v)\neq\emptyset$, while $\vec H_{\rm out}(i)=\emptyset$ for $i\in\La_0\cup\La_1$.

The following crucial theorem, proved in Subsection~\ref{S:contour}, links the maximal trajectory to Toom contours. In its essence, this goes back to part~3 of the proof of \cite[Thm~1]{Too80}, but we have reformulated things to a point where, putting the two texts besides each other, it is hard at first sight to spot the similarity.

\bt[Presence of a Toom contour]
Let\label{T:contour} $(\La,\Hi)$ be a typed dependence graph with $\sig\geq 1$ types of edges, let $\phhb$ be its associated monotone cellular automaton, and let $\ov x$ be its maximal trajectory. If $\ov x(i)=0$ for some $i\in\La$, then a Toom contour $(v_\circ,\Vi,\Ei,\psi)$ rooted at $i$ is present in $(\La,\Hi)$.
\et

We note that the converse of Theorem~\ref{T:contour} does not hold, i.e., the presence in $(\La,\Hi)$ of a Toom contour $(v_\circ,\Vi,\Ei,\psi)$ does not imply that $\ov x(i)=0$. This can be seen from Figure~\ref{fig:minexpl1}. In this example, if there would be no other defective sites apart from the sinks of the Toom contour, then the origin would have the value one. This is a difference with the Peierls arguments used in percolation theory, where the presence of a contour is a necessary and sufficient condition for the absence of percolation.

Theorem~\ref{T:contour} is one of our main results. In its essence, including the proof we have given, it goes back to the proof of \cite[Thm~1]{Too80}. In particular, in part~3 of that proof, by induction two ``trusses'', a ``set of classes'' and a non-oriented graph are constructed that have a number of properties numbered 1--7. These do not quite amount to our Theorem~\ref{T:contour}, but they list enough properties of these objects to allow Toom to prove his stability theorem. For example, property~6 in Toom's list states that the number of elements of the ``set of classes'' equals the number of ``forks'' plus one. These ``forks'' correspond to what we call sources other than the root, while the elements of the ``set of classes'' correspond to the sinks. Toom's property~6 therefore corresponds to the fact that the number of sinks equals the number of sources (including the root). A difference between our Toom contours and Toom's ``trusses'' is that he removes the sources other than the root from his graphs and connects the vertices adjacent to them by collection of directed edges, which he calls a ``fork''. Another important difference is that Toom reverses the direction of some arrows and doubles some other arrows with the aim of making sure that the same number of arrows enter and exit \emph{each} vertex, including the source at the root and the sinks. So in Toom's formulation, there no sources or sinks as the incoming and outgoing charges are equal at each vertex. This helps Toom prepare for the task of counting the number of contours of a certain size, but somewhat obscures why his property~6 holds. Instead of proving property~6 by observing that the number of sources equals the number of sinks, as we do, Toom obtains this property as a side result of his inductive construction of contours, which corresponds roughly to our proof of Theorem~\ref{T:contour}.

In spite of the significant changes in the formulation, our proof of Theorem~\ref{T:contour} still largely follows the original argument of Toom \cite{Too80}. A difference is our systematic use of the concept of minimal explanations, see Subsection~\ref{S:minexpl} below. The connection between Toom contours and minimal explanations has been observed before in \cite{Gac95,Pre07} but some remarks by these authors seem to suggest that the relation is more straightforward than it really is, a point we have sought to clarify in Theorem~\ref{T:embed} below. Toom cycles, which we will discuss in the next subsection, are not treated in \cite{Too80}.

\subsection{Toom cycles}\label{S:cyc}

We will give two proofs of Theorem~\ref{T:contour}: one that works for general $\sig\geq 1$, and another that works only for $\sig=2$, but that in this case gives some extra information that can sometimes be used to get sharper bounds. As Figure~\ref{fig:minexpl1} shows, Toom contours with two charges are essentially cycles. In the present subsection, we define Toom cycles, which are Toom contours with two charges that have some useful additional properties, and we formulate a theorem about the presence of Toom cycles in monotone cellular automata.

For $n\geq 2$, let $[n]:=\{0,\ldots,n-1\}$ equipped with addition modulo $n$. We define a \emph{cycle} of \emph{length} $n\geq 2$ to be an undirected graph $(V,E)$ with vertex set $V=[n]$ and edge set $E:=\big\{\{v,v+1\}:v\in[n]\big\}$. Similarly, we define a cycle of length $1$ to be the undirected graph $(V,E):=(\{0\},\emptyset)$. We define an \emph{oriented cycle} of \emph{length} $n\geq 1$ to be a directed graph $(V,\vec E)$ whose associated undirected graph $(V,E)$ is a cycle of length $n$ such that for each undirected edge $\{v,w\}\in E$, precisely one of the directed edges $(v,w)$ and $(w,v)$ is an element of $\vec E$ (but not both). In other words, this is a cycle in which each undirected edge has been given an orientation.

Each oriented cycle of length $n\geq 2$ naturally gives rise to a connected Toom graph $(\Vi,\Ei)$ with two charges by setting
\bc\label{VVV}
\dis V_\circ&:=&\dis\big\{v\in[n]:(v,v-1),(v,v+1)\in\vec E\big\},\\[5pt]
\dis V_\ast&:=&\dis\big\{v\in[n]:(v-1,v),(v+1,v)\in\vec E\big\},\\[5pt]
\dis V_1&:=&\dis\big\{v\in[n]:(v-1,v),(v,v+1)\in\vec E\big\},\\[5pt]
\dis V_2&:=&\dis\big\{v\in[n]:(v+1,v),(v,v-1)\in\vec E\big\}
\ec
and
\be
\vec E_1:=\big\{(v,w)\in\vec E:w=v+1\big\}
\quand
\vec E_2:=\big\{(v,w)\in\vec E:w=v-1\big\}.
\ee
Similarly, we may associate the oriented cycle of length one with the trivial Toom graph $(\Vi,\Ei)$ with two charges defined as $V_\circ=V_\ast:=\{0\}$ and $V_1=V_2=\vec E_1=\vec E_2:=\emptyset$. If $0\in V_\circ$, then we can take $v_\circ:=0$ to be the root. In view of this, connected rooted Toom graphs with two charges correspond (up to isomorphism) precisely to oriented cycles $(V,\vec E)$ of length $n\geq 1$ for which $0\in V_\circ$.

It is sometimes convenient to add the element $n$ to $V$ and to replace the oriented edge $(0,n-1)\in\vec E_2$ by $(n,n-1)$. Thus, we may identify an oriented cycle $(V,\vec E)$ of length $n\geq 1$ for which $0\in V_\circ$ with an oriented path of length $n$ for which $(0,1)\in\vec E_1$ and $(n,n-1)\in\vec E_2$. Similar to what we did in (\ref{toomroot}), it will be convenient to define
\be\label{toomroot2}
V'_\circ:=V_\circ\beh\{0\},\quad V'_1:=V_1\cup\{0\}\quand V'_2:=V_1\cup\{n\}.
\ee
In condition (ii) of the following definition, we equip $\{0,\ldots,n\}$ with the natural total order and we equip $\{1,\circ,2\}$ with the total order $1<\circ<2$.

\begin{defi}
Let\label{D:cycle} $\La$ be a countable set. A \emph{Toom cycle} in $\La$ is a triple $(V,\vec E,\psi)$ where $(V,\vec E)$ is an oriented cycle of length $n\geq 1$ such that $0\in V_\circ$, and $\psi:\{0,\ldots,n\}\to\La$\footnote{For ease of notation we write $\psi_v$ for $\psi(v)$ in case of Toom cycles.} is a map such that $\psi_0=\psi_n$ and
\begin{enumerate}
\item $\psi_v\neq\psi_w$ for each $v\in V_\ast$ and $w\in V$ with $v\neq w$,
\item if $\psi_v=\psi_w$ for some $v\in V'_s$ and $w\in V'_t$ with $s,t\in\{1,\circ,2\}$ and $s\leq t$, then $v\leq w$.
\end{enumerate}
We say that the Toom cycle $(V,\vec E,\psi)$ is \emph{rooted} at $\psi_0$.
\end{defi}

If $(V,\vec E,\psi)$ is a Toom cycle and $(v_\circ,\Vi,\Ei)$ is its associated rooted Toom graph, then $(v_\circ,\Vi,\Ei,\psi)$ is a Toom contour with two charges. We call this the Toom contour \emph{associated with} the Toom cycle $(V,\vec E,\psi)$. Applying property~(ii) with $s=t$ implies
\begin{itemize}
\item if $\psi_v=\psi_w$ for some $v,w\in V'_s$ with $s\in\{1,\circ,2\}$, then $v=w$,
\end{itemize}
so property~(ii) of Definition~\ref{D:cycle} implies property~(ii) of Definition~\ref{D:embed}. It is easy to see that it is strictly stronger, so not every Toom contour with two charges comes from a Toom cycle.

We next define what it means for a Toom cycle to be present in a typed dependence graph $(\La,\Hi)$.

\begin{defi}
Let\label{D:cycpres} $(\La,\Hi)$ be a typed dependence graph with 2 types of edges as in Definition~\ref{D:typdep}. We say that a Toom cycle $(V,\vec E,\psi)$ of length $n\geq 2$ is present in $(\La,\Hi)$ if:
\begin{enumerate}
\item $\dis\psi_v\in\La_0$ for all $v\in V_\ast$,
\item $\dis(\psi_v,\psi_w)\in\vec H_s$ for all $(v,w)\in\vec E_s$ with $v\in V'_s$ $(s=1,2)$,
\item $\dis(\psi_v,\psi_w)\in\vec H_{3-s}$ for all $(v,w)\in\vec E_s$ with $v\in V'_\circ$ $(s=1,2)$.
\end{enumerate}
For Toom cycles of length 1, only condition~(i) applies.
\end{defi}

Note that in condition~(iii) above, $3-s$ is $1$ if $s=2$ and $2$ if $s=1$, so this condition says that directed edges coming out of a source other than the root must use a directed edge of $(\La,\Hi)$ of the opposite charge. This condition is stronger than condition~(iii) of Definition~\ref{D:present}. One can check that our definition implies that if a Toom cycle is present in $(\La,\Hi)$, then its associated Toom contour is present in $(\La,\Hi)$ in the sense of Definition~\ref{D:present}, but because of our previous remark, the converse implication does not hold. One can check that the Toom contour with two charges in Figures \ref{fig:minexpl1} and \ref{fig:minexpl2} comes from a Toom cycle that is present in the strong sense of Definition~\ref{D:cycpres}.

In the same way as in (\ref{noLa1}), one can see that Definition~\ref{D:cycpres} implies that $\psi_v\not\in\La_1$ for all $v\in V$. In view of our previous remarks, the following theorem strengthens Theorem~\ref{T:contour} in the special case of two charges. Our proof of Theorem~\ref{T:cycle} (in Subsection~\ref{S:cycle}) will largely be independent of the proof of Theorem~\ref{T:contour}.

\bt[Presence of a Toom cycle]
Let\label{T:cycle} $(\La,\Hi)$ be a typed dependence graph with $2$ types of edges, let $\phhb$ be its associated monotone cellular automaton, and let $\ov x$ be its maximal trajectory. If $\ov x(i)=0$ for some $i\in\La$, then a Toom cycle rooted at $i$ is present in $(\La,\Hi)$.
\et

\subsection{A Peierls bound}\label{S:Peiabs}

Theorems \ref{T:contour} and \ref{T:cycle} can be used to prove upper bounds on the probability that the maximal trajectory of a random monotone cellular automaton takes the value zero in a given point. For concreteness, we formulate this as a theorem.

\bt[Peierls bound]
Let\label{T:ToomPei} $\Phi=(\Phi_i)_{i\in\La}$ be a random monotone cellular automaton and let $\ov X$ be its maximal trajectory. Let
\be
\La_r:=\big\{i\in\La:\Phi_i=\phh^r\big\}\quad(r=0,1)
\quand\La_\bullet:=\La\beh(\La_0\cup\La_1).
\ee
Let $\sig\geq 1$ be an integer and for each $i\in\La_\bullet$ and $\lis$, let ${\rm A}_{s,i}\in\Oi(\Phi_i)$. Define $\Hi=(\vec H_1,\ldots,\vec H_\sig)$ by
\be\label{Hs}
\vec H_s:=\big\{(i,j):i\in\La_\bullet,\ j\in {\rm A}_{s,i}\big\}\qquad(\lis).
\ee
Fix $i\in\La$ and let $\Ti_i$ denote the set of Toom contours rooted at $i$ (up to equivalence). Then
\be\label{ToomPei}
\P\big[\ov X(i)=0\big]\leq\sum_{T\in\Ti_i}\P\big[T\mbox{ is present in }(\La,\Hi)\big].
\ee
If $\sig=2$, then (\ref{ToomPei}) remains true if we restrict the sum to Toom cycles rooted at~$i$.
\et

\bpro
Let $\Psi$ be the random monotone cellular automaton associated with the random typed dependence graph $(\La,\Hi)$. Then in view of (\ref{phhFi}) we have $\Psi_i=\Phi_i$ for $i\in\La_0\cup\La_1$ while
\be
\Psi_i(x)=\bigvee_{s=1}^\sig\bigwedge_{j\in {\rm A}_{s,i}}x(j)\leq\bigvee_{A\in\Oi(\Phi_i)}\bigwedge_{j\in A}x(j)=\Phi_i(x)\qquad\big(i\in\La_\bullet,\ x\in\{0,1\}^\La\big).
\ee
Using this, it is easy to check (see Lemma~\ref{L:upcomp} below) that the maximal trajectories $\ov X$ of $\Phi$ and $\ov Y$ of $\Psi$ are ordered as $\ov Y\leq\ov X$ (pointwise). In particular,
\be\label{XleqY}
\P\big[\ov X(i)=0\big]\leq\P\big[\ov Y(i)=0\big].
\ee
By Theorems \ref{T:contour} and \ref{T:cycle}, the right-hand side of (\ref{XleqY}) can be bounded from above by the probability that there is a Toom contour or cycle present in $(\La,\Hi)$, which in turn can be estimated from above by the expected number of Toom contours or cycles.
\epro

\part{Proofs}

\section{Preliminaries}\label{S:prelim}

\subsection{Eroders}\label{S:erod}

In this subsection we prove Lemmas \ref{L:edge} and \ref{L:erode}.\med

\bpro[of Lemma~\ref{L:edge}]
It suffices to prove the claim for $t=1$. Fix $j\in\Z^d$ and set $j+A:=\{j+i:i\in A\}$ $(A\in\Oi(\phi))$. Then one has $\Psi_{0,t}(H^\ell_r)(j)=1$ if and only if there exists an $A\in\Oi(\phi)$ such that $\ell(k)\geq r$ for all $k\in j+A$. Equivalently, this says that
\be
\sup_{A\in\Oi(\phi)}\inf_{k\in j+A}\ell(k)\geq r.
\ee
Using (\ref{edge}) and linearity, we can rewrite this as $\ell(j)+\eps_\phi(\ell)\geq r$, which is equivalent to $j\in H^\ell_{r-\eps_\phi(\ell)}$.
\epro

We next prove Lemma~\ref{L:erode}. Our proof depends on the equivalence of (\ref{erosion}) and the eroder property, which is proved in \cite[Thm~1]{Pon13}. We recall that the fact that (\ref{erosion}) implies the eroder property has already been demonstrated below Lemma~\ref{L:erode}, so we depend on \cite[Thm~1]{Pon13} only for the converse implication.\med

\bpro[of Lemma~\ref{L:erode}]
In \cite[Lemma~12]{Pon13} it is shown\footnote{Since Ponselet discusses stability of the all-zero fixed point while we discuss stability of the all-one fixed point, in \cite{Pon13} the roles of zeros and ones are reversed compared to our conventions.} that (\ref{erosion}) is equivalent to the existence of a linear polar function $L$ of dimension $2\leq\sig\leq d+1$ and constants $\eps_1,\ldots,\eps_\sig$ such that $\sum_{s=1}^\sig \eps_s>0$ and for each $\lis$, there exists an $A_s\in\Oi(\phi)$ such that $\eps_s-L_s(i)\leq 0$ for all $i\in A_s$. It follows that
\be
\sum_{s=1}^\sig\sup_{A\in\Oi(\phi)}\inf_{i\in A}L_s(i)
\geq\sum_{s=1}^\sig\inf_{i\in A_s}L_s(i)\geq\sum_{s=1}^\sig \eps_s>0,
\ee
which shows that (\ref{erode}) holds. Assume, conversely, that (\ref{erode}) holds. Since $\Oi(\phi)$ is finite, for each $\lis$ we can choose $A_s\in\Oi(\phi)$ such that
\be\label{epss}
\eps_s:=\inf_{i\in A_s}L_s(i)=\sup_{A\in\Oi(\phi)}\inf_{i\in A}L_s(i).
\ee
Then (\ref{erode}) says that $\sum_{s=1}^\sig\eps_s>0$. Let $H_s:=\{z\in\R^d:L_s(z)\geq\eps_s\}$. By the definition of a linear polar function, $\sum_{s=1}^\sig L_s(z)=0$ for each $z\in\R^d$, and hence the condition $\sum_{s=1}^\sig\eps_s>0$ implies that for each $z\in\R^d$, there exists an $\lis$ such that $L_s(z)<\eps_s$. In other words, this says that $\bigcap_{s=1}^\sig H_s=\emptyset$. For each $\lis$, the set $A_s$ is contained in the half-space $H_s$ and hence the same is true for ${\rm Conv}(A_s)$, so we conclude that
\be
\bigcap_{s=1}^\sig{\rm Conv}(A_s)=\emptyset,
\ee
from which (\ref{erosion}) follows.
\epro

\subsection{The maximal trajectory}\label{S:maxtraj}

In this subsections, we prove Lemmas \ref{L:maxtraj} and \ref{L:numax}, as well as Lemma~\ref{L:upcomp} that has already been used in the proof of Theorem~\ref{T:ToomPei}.\med

\bpro[of Lemma~\ref{L:maxtraj}]
By symmetry, it suffices to show that there exists a trajectory $\ov x$ that is uniquely characterised by the property that each trajectory $x$ of $\phhb$ satisfies $x\leq\ov x$. Let $\La_n\sub\La$ be finite sets increasing to $\La$ and for each $n$, let $\phhb^n$ denote the monotone cellular automaton defined by
\be\label{phhcut}
\phh^n_i:=\bcase
\phh^1&i\in\La\beh\La_n\\
\phh_i&i\in\La_n,
\ecase
\ee
where $\phh^1$, defined in (\ref{phh01}), denotes the map that is constantly one. Since $\La_n$ is finite, it is easy to see that $\phhb^n$ has a unique trajectory $x^n$, which satisfies $x^n(i)=1$ for all $i\in\La\beh\La_n$. One has $x^n\geq x^{n+1}$ (coordinatewise) for each $n$ so the monotone limit $\ov x(i):=\lim_{n\to\infty}x^n(i)$ $(i\in\La)$ exists. It is straightforward to check that $\ov x$ is a trajectory of $\phhb$. If $x$ is any other trajectory of $\phhb$, then $x\leq x^n$ for all $n$ and hence $x\leq\ov x$.
\epro

\bpro[of Lemma~\ref{L:numax}]
By symmetry, it suffices to prove the claim for the upper invariant law. For each $n\geq 0$, let $\Phi^{n,p}$ denote the modified cellular automaton defined by
\be
\Phi^{n,p}_{i,t}:=\bcase
\phi^1&t\leq-n\\
\Phi^p_{i,t}&t>-n.
\ecase
\ee
Then it is easy to see that $\Phi^{n,p}$ has a unique trajectory $X^{n,p}$, which satisfies $X^{n,p}(i,t)=1$ for all $i\in\Z^d$ and $t\leq-n$. Exactly the same argument as in the proof of Lemma~\ref{L:maxtraj} shows that $X^{n,p}\to\ov X^p$ (pointwise) almost surely. The claim now follows from the observation that $(X^{n,p}(i,0))_{i\in\Z^d}$ is equally distributed with the random variable $X^p_n$ in the second formula of (\ref{toupper}).
\epro

\bl[Comparison of maximal trajectories]
Let\label{L:upcomp} $\phhb=(\phh_i)_{i\in\La}$ and $\psib=(\psi_i)_{i\in\La}$ be monotone cellular automata and let $\ov x$ and $\ov y$ denote their respective maximal trajectories. Assume that
\be\label{upcomp}
\phh_i(x)\leq\psi_i(x)\qquad\big(i\in\La,\ x\in\{0,1\}^\La\big).
\ee
Then $\ov x(i)\leq\ov y(i)$ $(i\in\La)$.
\el

\bpro
Define $\phhb^n$ and $\psib^n$ as in (\ref{phhcut}) and let $x^n$ and $y^n$ denote their unique trajectories. Then by induction, (\ref{upcomp}) implies that $x^n\leq y^n$ (pointwise), so taking the limit we obtain that $\ov x\leq\ov y$.
\epro

\subsection{Complete instability}\label{S:instab}

Let $\Phi^p$ be defined as in (\ref{Phip}) with $m=1$, i.e., $\Phi^p$ is a random perturbation of the deterministic cellular automaton $\Phi^0$ that applies the same nonconstant local monotone map $\phi_1=\phi$ in each space-time point. Let $\ov\rho(p)$, defined in (\ref{ovrho}), denote the density of its upper invariant law. Our Theorem~\ref{T:detspeed} implies as a special case the difficult part of Toom's stability theorem (Theorem~\ref{T:Toom}), which says that $\Phi^0$ is stable if $\phi$ is an eroder. In the present subsection, we complement this by proving the ``easy'' part of Toom's stability theorem, which says that $\Phi^0$ is completely unstable if $\phi$ is not an eroder.

\bl[Complete instability]
If\label{L:instab} $\phi$ is not an eroder, then $\ov\rho(p)=0$ for all $p>0$.
\el

\bpro
By translation invariance, it suffices to prove that for each $p>0$, the Markov chain $(X^p_t)_{t\geq 0}$ defined in (\ref{Markov}) and started in the initial state $X^p_0=\un 1$ satisfies
\be\label{tozero}
\P^{\un 1}\big[X^p_t(0)=1\big]\asto{t}0.
\ee
Since $\phi$ is not an eroder, there exists configuration $x\in\{0,1\}^{\Z^d}$ containing finitely many zeros such that $\Psi^t_\phi(x)\neq\un 1$ for all $t\geq 0$. This allows us to choose for each $t\geq 0$ a point $i_t\in\Z^d$ such that $\Psi^t_\phi(x)(i_t)=0$. Let us write $x(i)=1-1_A(i)$ $(i\in\Z^d)$ where $A\sub\Z^d$ is a finite set. Let $A-i_t:=\{i-i_t:i\in A\}$ $(t\geq 0)$. Then monotonicity implies that
\be
X^p_t(0)=0\quad\mbox{a.s.\ on the event that}\quad
\exists 0<s\leq t\mbox{ s.t.\ }\Phi^p_{s,i}=\phi_0\ \forall i\in A-i_s.
\ee
It follows that
\be
\P^{\un 1}\big[X^p_t(0)=1\big]\leq(1-p^{|A|})^t\qquad(t\geq 0),
\ee
which proves (\ref{tozero}).
\epro

\section{Construction of Toom contours}\label{S:Toomconstr}

\subsection{Minimal explanations}\label{S:minexpl}

This section is devoted to the proofs of Theorems \ref{T:contour} and \ref{T:cycle}, which can be found in Subsections \ref{S:contour} and \ref{S:cycle} below. In the present subsection, we prepare for these proofs by giving a formal definition of the minimal explanations that have already been mentioned several times, and investigating their properties.

\emph{Minimal zero-sets} of a monotone local map $\phi$ are defined analogously to the minimal one-sets of (\ref{Aphi}), i.e., these are minimal elements of the set of all finite $Z\sub\Z^d$ with the property that $\phi(1-1_Z)=0$. We let $\Zi(\phi)$ denote the set of all minimal zero-sets of $\phi$. In analogy with (\ref{Aphi}), each monotone local map $\phh:\{0,1\}^\La\to\{0,1\}$ can be written as
\be\label{AZphh}
\phh(x)=\bigvee_{A\in\Oi(\phh)}\bigwedge_{i\in A}x(i)=\bigwedge_{Z\in\Zi(\phh)}\bigvee_{i\in Z}x(i)\qquad\big(x\in\{0,1\}^\La\big).
\ee
In particular, if $\phh^0$ and $\phh^1$ are the constant maps defined in (\ref{phh01}), then $\Zi(\phh^0)=\{\emptyset\}$ and $\Zi(\phh^1)=\emptyset$. For monotone local maps $\phh,\phh'$, we write
\be
\phh\leq\phh'\quad\desd\quad\phh(x)\leq\phh'(x)\ \forall x\in\{0,1\}^\La
\quand
\phh\pre\phh'\quad\desd\quad\Zi(\phh)\supset\Zi(\phh').
\ee
It is easy to see that $\phh\pre\phh'$ implies $\phh\leq\phh'$, but not the other way around. For monotone cellular automata $\phhb=(\phh_i)_{i\in\La}$ and $\phhb'=(\phh'_i)_{i\in\La}$, we write $\phhb\leq\phhb'$ if and only if $\phh_i\leq\phh'_i$ for all $i\in\La$, and similarly, we write $\phhb\pre\phhb'$ if and only if $\phh_i\pre\phh'_i$ for all $i\in\La$.

\begin{defi}
Let\label{D:minexpl} $\phhb=(\phh_i)_{i\in\La}$ be a monotone cellular automaton and let $0\in\La$. By definition, a \emph{minimal explanation} for $0$ is a monotone cellular automaton $\phhb'$ such that:
\begin{enumerate}
\item $\phhb\pre\phhb'$ and the maximal trajectory $\ov x'$ of $\phhb'$ satisfies $\ov x'(0)=0$.
\item If a monotone cellular automaton $\phhb''$ satisfies $\phhb'\pre\phhb''$ and the maximal trajectory $\ov x''$ of $\phhb''$ satisfies $\ov x''(0)=0$, then $\phhb'=\phhb''$.
\end{enumerate}
\end{defi}

\bl[Minimal explanations]
Let\label{L:minexpl} $\phhb=(\phh_i)_{i\in\La}$ be a monotone cellular automaton and let $0\in\La$. Then the the maximal trajectory $\ov x$ of $\phhb$ satisfies $\ov x(0)=0$ if and only if there exists a minimal explanation $\phhb'$ for $0$.
\el

\bpro
If there exists a minimal explanation $\phhb'$ for $0$ and $\ov x$ and $\ov x'$ denote the maximal trajectories of $\phhb$ and $\phhb'$, respectively, then $\phhb\pre\phhb'$ implies $\phhb\leq\phhb'$ which implies $\ov x\leq\ov x'$ and hence in particular $\ov x(0)\leq\ov x'(0)=0$. This shows that $\ov x(0)=0$ if there exists a minimal explanation for $0$.

Assume, conversely, that $\ov x(0)=0$. Let $\La_n\sub\La$ be finite sets increasing to $\La$, let $\phhb^n$ denote the monotone cellular automata defined in (\ref{phhcut}), and let $x^n$ denote the unique trajectory of $\phhb^n$. We have seen in the proof of Lemma~\ref{L:maxtraj} that $\lim_{n\to\infty}x^n(0)=\ov x(0)$ so we can choose $n$ large enough such that $x^n(0)=0$. It is clear from the definition that $\phhb\pre\phhb^n$. We can now step by step replace $\phhb^n$ by larger monotone cellular automata with respect to the order $\pre$ as long as it is possible to do so without losing the property that the trajectory is zero in $0$. Since $\Zi(\phh^n_j)=\emptyset$ for all but finitely many $j$ and since $\Zi(\phh^n_j)$ is finite for each $j$, this process ends after a finite number of steps, leading to a minimal explanation for $0$.
\epro

Our next proposition describes the structure of minimal explanations. In point~(iii) below, we use the convention that the maximum over an empty set is zero. We call the finite directed graph $(U,\vec G)$ from Proposition~\ref{P:minexpl} the \emph{explanation graph} associated with the minimal explanation $\phhb'$. The picture on the right in Figure~\ref{fig:minexpl2} shows an example of such an explanation graph, or rather the undirected graph $(U,G)$ associated with $(U,\vec G)$.

\bp[Explanation graphs]
Let\label{P:minexpl} $\phhb=(\phh_i)_{i\in\La}$ be a monotone cellular automaton and let $(\La,\vec H)$ be its dependence graph, as defined in (\ref{vecF}). Let $0\in\La$ and let $\phhb'$ be a minimal explanation for $0$. Then there exists a finite subgraph $(U,\vec G)$ of $(\La,\vec H)$ with the following properties.
\begin{enumerate}
\item The maximal trajectory $\ov x'$ of $\phhb'$ satisfies $\ov x'(i)=0$ if and only if $i\in U$.
\item $\phh'_i=\phh^1$ if $i\not\in U$.
\item $\dis\phh'_i(x)=\!\!\bigvee_{j:\,(i,j)\in\vec G}\!\!x(j)\quad(x\in\{0,1\}^\La)$ if $i\in U$.
\item For each $j\in U\beh\{0\}$, there exists an $i\in U$ such that $(i,j)\in\vec G$.
\end{enumerate}
For $i\in U$, the following statements are equivalent: 1.\ $\phh_i=\phh^0$, 2.\ $\phh'_i=\phh^0$, 3.\ $\vec G_{\rm out}(i)=\emptyset$.
\ep

\bpro
Define $U:=\{j\in\La:\phh'_j\neq\phh^1\}$. We claim that $U$ is finite. This follows from the argument we have already seen in the proof of Lemma~\ref{L:minexpl}: if $\La_n\sub\La$ are finite sets increasing to $\La$, then for large enough $n$ we can replace $\phh'_j$ by $\phh^1$ for all $j\not\in\La_n$ without affecting the fact that the maximal trajectory is zero in $0$. By the maximality property of $\phhb'$, this then implies that $\phh'_j=\phh^1$ for all $j\not\in\La_n$.

It is clear from our definition of $U$ that $\ov x'(j)=1$ for all $j\in\La\beh U$. On the other hand, we cannot have $\ov x'(j)=1$ for some $j\in U$, since in that case we could replace $\phh'_j$ by $\phh^1$ while preserving the fact that the maximal trajectory is zero in $0$, which contradicts the maximality property of $\phhb'$. This proves property~(i). Property~(ii) is immediate from our definition of $U$.

We claim that for each $i\in U$, there exists a finite set $Z_i\sub U$ such that $\Zi(\phh'_i)=\{Z_i\}$. Indeed, property~(i) and (\ref{AZphh}) imply that for each $i\in U$ there exists a $Z\in\Zi(\phh'_i)$ such that $\ov x'(k)=0$ for all $k\in Z$, which by (i) implies $Z\sub U$. If $\Zi(\phh'_i)$ contains other elements apart from $Z$, then we can throw these away while preserving the fact that the maximal trajectory is zero in $0$, contradicting the maximality property of $\phhb'$. Now setting
\be
\vec G:=\big\{(i,j):i\in U,\ j\in Z_i\big\}
\ee
defines a set of directed edges such that $\vec G\sub\vec H$ and property~(iii) holds. Note that in line with earlier conventions, we allow for the case that $Z_i=\emptyset$ and $\phh_i=\phh^0$.

Property~(iv) follows from the fact that if $j\in U\beh\{0\}$ and there exists no $i\in U$ such that $j\in Z_i$, then by property~(ii) and (\ref{AZphh}) we can replace $\phh'_j$ by $\phh^1$ while preserving the fact that the maximal trajectory is zero in all points of $U\beh\{j\}$, contradicting the maximality property of $\phhb'$.

To prove the final statement of the proposition, we observe that if $\phh_i=\phh^0$, then $\vec G_{\rm out}(i)\sub\vec H_{\rm out}(i)=\emptyset$, so 1.\ implies 3. By property~(iii), 3.\ implies 2., which by the fact that $\phhb\prec\phhb'$ in turn implies 1.
\epro

In the special case that the monotone cellular automaton $\phhb$ is defined in terms of a typed dependence graph $(\La,\Hi)$ as in Definition~\ref{D:typdep}, we can strengthen Proposition~\ref{P:minexpl} as follows. We call the typed directed graph $(U,\Gi)$ from the following proposition a \emph{typed explanation graph} associated with the minimal explanation $\phhb'$. In general,  $(U,\Gi)$ is not uniquely determined by $\phhb'$.

\bp[Typed explanation graphs]
Let\label{P:typexpl} $(\La,\Hi)$ be a typed dependence graph with $\sig\geq 1$ types of edges and let $\phhb$ be its associated monotone cellular automaton. Let $0\in\La$ and let $\phhb'$ be a minimal explanation for $0$. Then there exists a finite typed subgraph $(U,\Gi)$ of $(\La,\Hi)$ such that:
\begin{enumerate}
\item $\phh'_i=\phh^1$ if $i\not\in U$,
\item $\phh'_i=\phh^0$ if $i\in U_\ast:=\{i\in U:\phh_i=\phh^0\}$,
\item for each $i\in U\beh U_\ast$ and $\lis$, there exists a $j_s(i)\in U$ such that $\vec G_{s,{\rm out}}(i)=\{j_s(i)\}$,
\item $\dis\phh'_i(x)=\bigvee_{s=1}^\sig x\big(j_s(i)\big)$ $(x\in\{0,1\}^\La)$ if $i\in U\beh U_\ast$.
\end{enumerate}
The untyped directed graph $(U,\vec G)$ associated with $(U,\Gi)$ is the explanation graph associated with the minimal explanation $\phhb'$.
\ep

\bpro
Let $(U,\vec G)$ be the explanation graph associated with the minimal explanation $\phhb'$ and for each $i\in U$, let $Z_i:=\{j\in U:(i,j)\in\vec G\}$. Then property~(iii) of Proposition~\ref{P:minexpl} says that
\be
\phh'_i(x)=\bigvee_{j\in Z_i}x(j)\qquad\big(i\in U,\ x\in\{0,1\}^\La\big),
\ee
so $\Zi(\phh'_i)=\{Z_i\}$ and hence $Z_i\in\Zi(\phh_i)$ by the fact that $\phh_i\pre\phh'_i$. Let
\be\label{Asi}
A^s_i:=\big\{j\in\La:(i,j)\in\vec H_s\big\}\qquad(\lis,\ i\in\La).
\ee
Recall from Definition~\ref{D:typdep} that $\La_\bullet=\big\{i\in\La:\phh_i\not\in\{\phh^0,\phh^1\}\big\}$ and that
\be\label{phU}
\phh_i(x)=\bigvee_{s=1}^\sig\bigwedge_{j\in A^s_i}x(j)\qquad\big(i\in U,\ x\in\{0,1\}^\La\big).
\ee
We claim that
\be\label{AcapZ}
A^s_i\cap Z\neq\emptyset\qquad\big(Z\in\Zi(\phh_i),\ \lis,\ i\in\La_\bullet\big).
\ee
Indeed, if we would have $A^s_i\cap Z=\emptyset$ for some $Z\in\Zi(\phh_i)$, $\lis$, and $i\in\La_\bullet$, then $1=\phi_i(1_{A^s_i})\leq\phi_i(1-1_Z)=0$, which is a contradiction. By (\ref{AcapZ}), for each $\lis$ and $i\in\La_\bullet$, we can choose $j_s(i)\in A^s_i\cap Z_i$. Let $Z'_i:=\{j_s(i):\lis\}$. Then clearly $Z'_i\sub Z_i$. We claim that in fact $Z'_i=Z_i$. Indeed, since $j_s(i)\in A^s_i$ $(\lis)$, formula (\ref{phU}) shows that $\phh_i(1-1_{Z'_i})=0$. Since $Z'_i\sub Z_i$, by the minimality of the latter, we conclude that $Z'_i=Z_i$. As a result, defining a typed directed graph $(U,\Gi)$ with $\sig$ types of edges by setting
\be
\vec G_s:=\big\{\big(i,j_s(i)\big):i\in U\cap\La_\bullet\big\}\qquad(\lis),
\ee
we have that $(U,\vec G)$ is the untyped directed graph associated with $(U,\Gi)$ and
\be
\phh'_i(x)=\bigvee_{j\in Z'_i}x(j)=\bigvee_{s=1}^\sig x\big(j_s(i)\big)
\qquad\big(i\in U\cap\La_\bullet,\ x\in\{0,1\}^\La\big).
\ee
Now properties (i)--(iv) follow from properties (ii) and (iii) of Proposition~\ref{P:typexpl}, while (\ref{Asi}) shows that $(U,\Gi)$ is a typed subgraph of $(\La,\Hi)$.
\epro

\subsection{Toom contours}\label{S:contour}

In this subsection, we prove Theorem~\ref{T:contour}. We fix a typed dependence graph $(\La,\Hi)$ with $\sig\geq 1$ types of edges. We let $\phhb$ denote its associated monotone cellular automaton and let $\ov x$ denote its maximal trajectory. We fix an element $0\in\La$ and assume that $\ov x(0)=0$. We need to prove the presence in $(\La,\Hi)$ of a Toom contour $(v_\circ,\Vi,\Ei,\psi)$ rooted at $0$. Since $\ov x(0)=0$, by Lemma~\ref{L:minexpl}, there exists a minimal explanation $\phhb'$ for $0$, and by Proposition~\ref{P:typexpl}, there exists a typed explanation graph $(U,\Gi)$ associated with $\phhb'$. We will derive Theorem~\ref{T:contour} from the following theorem. Recall Definition~\ref{D:embed} of an embedding of a rooted Toom graph.

\begin{figure}[htb]
\begin{center}
\includegraphics[width=\textwidth]{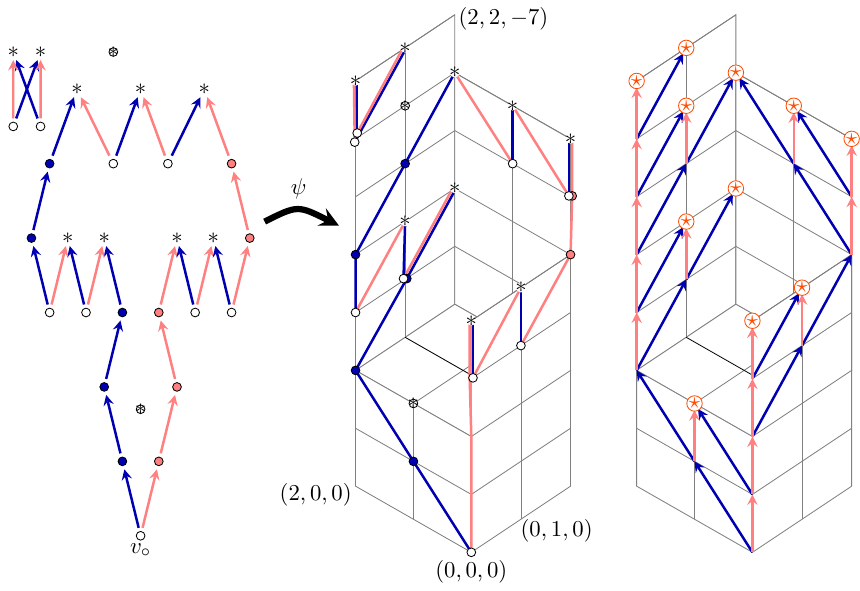}
\caption{Embedding of a rooted Toom graph inside a typed explanation graph. On the right: a typed explanation graph $(U,\Gi)$ associated with a minimal explanation for $(0,0,0)$ in the sense of Proposition~\ref{P:typexpl}. On the left and in the middle: embedding of a rooted Toom graph in $(U,\Gi)$ in the sense of Theorem~\ref{T:embed}. The connected component of this Toom graph containing the root is a Toom contour rooted at $(0,0,0)$ (compare Figure~\ref{fig:minexpl1}).}
\label{fig:minexpl2}
\end{center}
\end{figure}

\bt[Toom graph embedded in explanation graph]
Let\label{T:embed} $(U,\Gi)$ be a typed explanation graph associated with a minimal explanation $\phhb'$ for $0$. Then there exists a rooted Toom graph $(v_\circ,\Vi,\Ei)$ and an embedding $\psi$ of $(v_\circ,\Vi,\Ei)$ in $U$ such that $\psi(v_\circ)=0$ and
\begin{enumerate}
\item $\dis\psi(V_\ast)=U_\ast$,
\item $\dis\big(\psi(v),\psi(w)\big)\in\vec G_s$ for all $(v,w)\in\vec E^\bullet_s$ $(\lis$),
\item $\dis\big(\psi(v),\psi(w)\big)\in\vec G$ for all $(v,w)\in\vec E^\circ$,
\end{enumerate}
where $\vec E^\bullet_s$ and $\vec E^\circ$ are defined in (\ref{Ecirc}).
\et

To see that Theorem~\ref{T:embed} implies Theorem~\ref{T:contour}, it suffices to observe that if $(V',\Ei')$ is the connected component containing $v_\circ$ of the Toom graph $(\Vi,\Ei)$ from Theorem~\ref{T:embed}, and $\psi'$ is the restriction of $\psi$ to $V'$, then $(v_\circ,V',\Ei',\psi')$ is a Toom contour rooted at $0$ that is present in $(\La,\Hi)$. Note that when we restrict ourselves to the connected component containing the root, property~(i) of Theorem~\ref{T:embed} must be weakened to $\psi(V_\ast)\sub U_\ast$, which is all that is needed to satisfy Definition~\ref{D:present}~(i). Theorem~\ref{T:embed} is demonstrated in Figure~\ref{fig:minexpl2}.

We observe that Theorem~\ref{T:embed} is trivial if $|U|=1$, since in this case $0$ is a defective site and we can take for $(\Vi,\Ei)$ the trivial Toom graph that consists of a single isolated vertex. We assume therefore from now on that $|U|\geq 2$. In this case, $0\not\in U_\ast$.

The proof of Theorem~\ref{T:embed} needs some preparations. In any directed graph $(V,\vec E)$, for two vertices $i,j\in V$, we write $i\leadsto j$ if there exists $i=i_0,\ldots,i_n=j$ such that $(i_{k-1},i_k)\in\vec E$ $(1\leq k\leq n)$. By definition, a \emph{time-ordering} of $(V,\vec E)$ is an enumeration $V=\{i_1,\ldots,i_N\}$ of its vertices such that for each $1\leq n\leq N$, there are no $k,l$ with $k<n\leq l$ and $(i_l,i_k)\in\vec E$. Note that since $(\La,\vec H)$ is acyclic, the same is true for $(U,\vec G)$.

\bl[Time-ordering]
Each\label{L:timord} finite acyclic directed graph has a time-ordering. For an explanation graph $(U,\vec G)$, we can choose a time-ordering such that $i_1=0$ and $U\beh U_\ast=\{i_1,\ldots,i_m\}$ for some $1\leq m\leq|U|$.
\el

\bpro
For any acyclic directed graph, the relation $\leadsto$ is a partial order on $V$; in particular, there cannot exist $i,j\in V$ with $i\neq j$ such that $i\leadsto j$ and $j\leadsto i$, since this would imply the existence of a cycle in $(V,\vec E)$. We can now inductively construct a time-ordering $i_1,i_2,\ldots$ by choosing for $i_n$ a minimal element of $V\beh\{i_1,\ldots,i_{n-1}\}$.

It follows from Proposition~\ref{P:minexpl}~(iv) that $0\leadsto i$ for all $i\in U$, so by the fact that $(U,\vec G)$ is acyclic we have $i\not\leadsto 0$ for all $i\in U\beh\{0\}$. Thus, $0$ is a minimal element of $U$ with respect to the partial order $\leadsto$ and we can construct the time-ordering starting with $i_1=0$. Since elements of $U_\ast$ have no outgoing edges, we can also first construct a time-ordering of $U\beh U_\ast$ and then add the elements of $U_\ast$ in any order.
\epro

From now on, we fix a typed explanation graph $(U,\Gi)$ associated with a minimal explanation $\phhb'$ for $0$, as well as a time-ordering of the associated untyped explanation graph $(U,\vec G)$ with the properties described in Lemma~\ref{L:timord}. We let $m:=|U\setminus U_\ast|$ and we adopt the following definitions.

\begin{defi}
For\label{D:timepart} each $1\leq n\leq m$, we set $U^-_n:=\{i_1,\ldots,i_n\}$ and $U^+_n:=U\beh U^-_n$. We call
\be
\pa U^-_n:=\big\{j\in U^+_n:\exists i\in U^-_n\mbox{s.t.\ }(i,j)\in\vec G\big\}
\ee
the \emph{boundary} of $U^-_n$. We equip $\pa U^-_n$ with the structure of an unoriented graph in which two elements $i,j\in\pa U^-_n$ are neighbours, denoted $i\approx j$, if there exists a $k\in U$ such that $i\leadsto k$ and $j\leadsto k$. We write $i\sim j$ if $i,j\in\pa U^-_n$ lie in the same connected component of this graph.
\end{defi}

The following lemma says that the number of connected components on the boundary $\pa U^-_n$ is non-decreasing in $n$. Note that at the end, when $n=m$, we have $\pa U^-_m=U_\ast$ and each element of $\pa U^-_m$ forms a connected component on its own. Therefore, starting from a single connected component at $n=1$ the boundary gradually breaks up into smaller and smaller connected components.

\bl[Break-up of boundary]
For\label{L:breakup} each $1<n\leq m$, if $C$ is a connected component of $\pa U^-_{n-1}$ and $i_n\not\in C$, then $C$ is also a connected component of $\pa U^-_n$. Each connected component of $\pa U^-_n$ that is not a connected component of $\pa U^-_{n-1}$ contains a vertex $j$ such that $(i_n,j)\in\vec G$.
\el

\bpro
We first prove that a connected component $C$ of $\pa U^-_{n-1}$ that does not contain $i_n$ is also a connected components of $\pa U^-_n$. For each $i,j\in C$, there exist $i(0),\ldots,i(k)\in C$ such that $i(0)\approx\cdots\approx i(k)$ with $i(0)=i$ and $i(k)=j$, which implies that $i\sim j$ in $U^-_n$. This shows that $C$ is contained in some connected component $C'$ of $\pa U^-_n$. We need to show that $C=C'$. Assume that conversely, $C'$ is strictly larger than $C$. Then we can find $i\in C$ and $j\in C'\beh C$ such that $i\approx j$. Since $j\in\pa U^-_n$ we must have either $j\in\pa U^-_{n-1}$ or $(i_n,j)\in\vec G$ (possibly both). If $j\in\pa U^-_{n-1}$, then $i\approx j$ implies $j\in C$ which contradicts our assumptions. However, if $(i_n,j)\in\vec G$, then $i\approx j$ implies $i\approx i_n$ which also contradicts our assumptions, since $C$ does not contain $i_n$.

To prove the second claim of the lemma, assume that $C$ is a connected component of $\pa U^-_n$ that is not a connected component of $\pa U^-_{n-1}$. Let $i$ be any element of $C$. If $(i_n,i)\in\vec G$ we are done. In the opposite case, $i\in\pa U^-_{n-1}$. Since $C$ is not a connected component of $\pa U^-_{n-1}$, by what we have already proved, $i$ must lie in the connected component of $\pa U^-_{n-1}$ that contains $i_n$, so there exist $i(0),\ldots,i(k)\in\pa U^-_{n-1}$ with $i(0)=i$, $i(k)=i_n$, and $i(0)\approx\cdots\approx i(k)$. Now there exists a $j'\in U$ such that $i(k-1)\leadsto j'$ and $i_n\leadsto j'$. If $j'\neq i_n$, then let $j$ be the first vertex after $i_n$ on the path from $i_n$ to $j'$, and if $j'=i_n$, then choose for $j$ any vertex with $(i_n,j)\in\vec G$. In either case, we then have $i(k-1)\approx j$ which implies that $j\in C$, and clearly $(i_n,j)\in\vec G$.
\epro

The idea of the proof of Theorem~\ref{T:embed} is to use a time-ordering of $(U,\vec G)$ as in Lemma~\ref{L:timord} to step by step build a Toom graph inside $(U,\vec G)$. In each step, the number of sources that have already been introduced is equal to the number of connected components of $\pa U^-_n$, and each connected component contains precisely one charge of each type. To formulate this idea precisely, we need one more definition.

\begin{defi}
For\label{D:ray} $\lis$, we define a \emph{spoke of charge} $s$ to be a sequence $\big(i(0),\ldots,i(k)\big)$ of vertices in $U$ such that $k\geq 1$, $i(k)\in U_\ast$, $\big(i(0),i(1)\big)\in\vec G$, and $\big(i(l-1),i(l)\big)\in\vec G_s$ for all $2\leq l\leq k$. We say that a spoke $\big(i(0),\ldots,i(k)\big)$ \emph{intersects} a set $V\sub U$ if $i(l)\in V$ for some $0\leq l\leq k$. A \emph{pole} at vertex $i\in U$ is a collection $\big(i_s(0),\ldots,i_s(k_s)\big)_{\lis}$ of spokes of charges $\lis$, respectively, such that $i_s(0)=i$ for all $\lis$.
\end{defi}

\bpro[of Theorem~\ref{T:embed}]
We have already shown that the statement is trivial if $0$ is a defective site, so we continue assuming that $|U|>1$ and $0\not\in U_\ast$. We fix a time-ordering of $(U,\vec G)$ as in Lemma~\ref{L:timord} and define $U^\pm_n$ as in Definition~\ref{D:timepart}. Let $N(n)$ denote the number of connected components of $\pa U^-_n$ $(1\leq n\leq m)$. It follows from Lemma~\ref{L:breakup} that $N(n)$ increases to $N(m)=|U_\ast|$. We will show by induction that for each $1\leq n\leq m$, it is possible to construct poles
\be
\big(i^r_s(0),\ldots,i^r_s(k^r_s)\big)_{\lis}\qquad\big(1\leq r\leq N(n)\big)
\ee
at vertices $i^1,\ldots,i^{N(n)}\in U^-_n$ such that
\begin{enumerate}
\item $i^1=0$ and $\big(i^1_s(0),i^1_s(1)\big)\in\vec G_s$ $(\lis)$,
\item for each connected component $C$ of $\pa U^-_n$ and for each $\lis$, there exists precisely one $1\leq r\leq N(n)$ such that the spoke $\big(i^r_s(0),\ldots,i^r_s(k^r_s)\big)$ intersects $C$.
\end{enumerate}
We start by proving the claim for $n=1$. By Proposition~\ref{P:typexpl}~(iii), for each $\lis$, at each $i\in U\beh U_\ast$ there is precisely one outgoing edge of charge $s$. Thus, for each $\lis$, there starts a unique spoke $\big(i^1_s(0),\ldots,i^1_s(k^1_s)\big)$ at $0$ such that $\big(i^1_s(l-1),i^1_s(l)\big)\in\vec G_s$ for all $1\leq l\leq k_s$, and these spokes together form a pole at $0$ such that (i) holds. If $\pa U^-_1$ has only one connected component, then (ii) also holds and we are now done. In the opposite case, we can add additional poles at 0 so that (ii) holds.

We now continue by induction on $n$. We will show that by adding poles, we can make sure (ii) remains valid as we increase $n$. Since (i) also obviously stays true if we add poles, this then completes the proof that (i) and (ii) can be satisfied for all $n$. Assume that we have poles at vertices $i^1,\ldots,i^{N(n-1)}\in U^-_{n-1}$ such that conditions (i) and (ii) are satisfied. By Lemma~\ref{L:breakup}, if $C$ is a connected component of $\pa U^-_{n-1}$ that does not contain $i_n$, then $C$ is also a connected component of $\pa U^-_n$, so for such a connected component $C$ condition~(ii) remains satisfied even without adding new poles. Let $C_1,\ldots,C_k$ be the other connected components of $\pa U^-_{n-1}$, which contain all vertices of the connected component of $\pa U^-_{n-1}$ that contains $i_n$, except $i_n$ itself, as well as all vertices $j\in U$ such that $(i_n,j)\in\vec G$. Using this and the induction hypothesis, we see that for each $\lis$, there exists precisely one $1\leq r\leq N(n-1)$ such that the spoke $\big(i^r_s(0),\ldots,i^r_s(k^r_s)\big)$ intersects $C_1\cup\cdots\cup C_k$. By Lemma~\ref{L:breakup}, each of the components $C_1,\ldots,C_k$ contains an element $j$ such that $(i_n,j)\in\vec G$. Using this and the fact that the edge between the first two vertices along each pole is in $\vec G$, we see that we can add additional poles in $i_n$ so that condition~(ii) is satisfied for $C_1,\ldots,C_k$. This completes the induction step.

In particular, setting $n=m$, we have now shown that it is possible to construct poles
\be
\big(i^r_s(0),\ldots,i^r_s(k^r_s)\big)_{\lis}\qquad\big(1\leq r\leq|U_\ast|\big)
\ee
at vertices in $U\beh U_\ast$ such that
\begin{enumerate}
\item $i^1=0$ and $\big(i^1_s(0),i^1_s(1)\big)\in\vec G_s$ $(\lis)$,
\item for each $i\in U_\ast$ and for each $\lis$, there exists precisely one $1\leq r\leq|U_\ast|$ such that the spoke $\big(i^r_s(0),\ldots,i^r_s(k^r_s)\big)$ ends in $i^r_s(k^r_s)=i$.
\end{enumerate}
It is straightforward to check that these poles together define a Toom graph $(\Vi,\Ei)$ that is embedded in $(U,\Gi)$ in such a way that conditions (i)--(iii) of the theorem are satisfied. Indeed, each pole corresponds to a source and its $\sig$ spokes to the $\sig$ charges emerging from the source. To see that $\psi$ satisfies condition~(ii) of Definition~\ref{D:embed} of an embedding of a rooted Toom graph, one uses the fact that if two spokes of the same charge would enter the same vertex, then these spokes would have to be equal starting from that vertex, which would lead to two spokes of the same charge ending in the same defective site, contradicting point~(ii) above. Also, since there are no incoming edges at $0$ in the explanation graph, we never add internal vertices that overlap with the root. The fact that conditions (ii) and (iii) of the theorem are satisfied follows our definition of a spoke of charge $s$, which includes the condition that $\big(i(l-1),i(l)\big)\in\vec G_s$ for all $2\leq l\leq k$, as well as point~(i) above.
\epro

\subsection{Toom cycles}\label{S:cycle}

\begin{figure}[htb]
\begin{center}
\includegraphics[width=\textwidth]{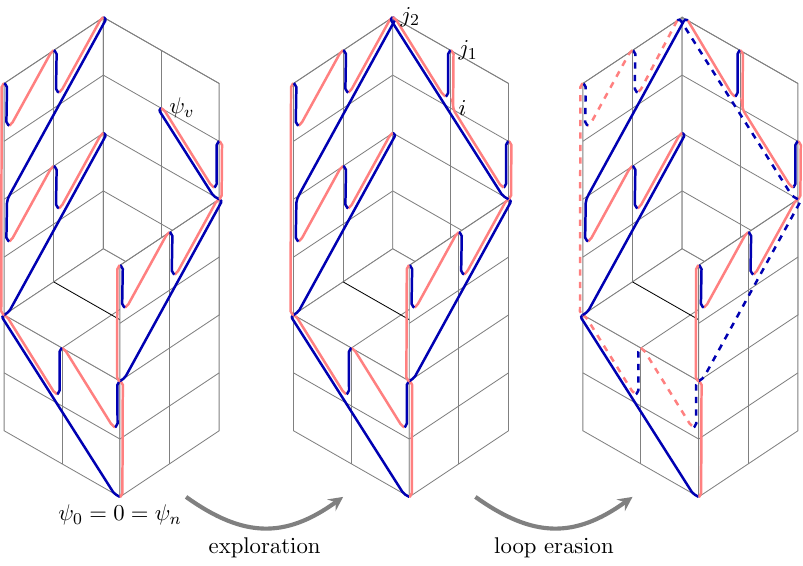}
\caption{The process of exploration and loop erasion. The Toom cycle is constructed on the explanation graph of Figure~\ref{fig:minexpl2}. We can see that in the Toom cycle on the left $v$ is a sink, but $\psi_v=i$ is not a defective site. In the exploration step, $v$ is replaced by two internal vertices, one of each charge, and two new sinks are added to the cycle at the positions  $j_1$ and $j_2$. This leads to the new sink at $j_2$ overlapping with a preexisting sink. In the loop erasion step, this is resolved by erasing the part of the cycle between the first and second visit to $j_2$.}
\label{fig:looperas}
\end{center}
\end{figure}

In this subsection, we prove Theorem~\ref{T:cycle}. As in the proof of Theorem~\ref{T:contour}, we will construct the Toom cycle inside a typed explanation graph. Apart from this similarity, the proof will be completely different. The proof is based on an inductive construction based on two steps, exploration and loop erasion, that are illustrated in Figure~\ref{fig:looperas}.\med

\bpro[of Theorem~\ref{T:cycle}]
We fix a typed explanation graph $(U,\Gi)$ as in Proposition~\ref{P:typexpl}, with two types of edges. By Proposition~\ref{P:minexpl}~(iv), for each $i\in U$, there exist $i_0,\ldots,i_n\in U$ with $i_0=0$ and $i_n=i$ such that $(i_{k-1},i_k)\in\vec G$ for all $1\leq k\leq n$. We let ${\rm dist}(i)$ denote the smallest integer $n$ for which such $i_0,\ldots,i_n$ can be found, i.e., ${\rm dist}(i)$ is the length of the shortest directed path in $(U,\vec G)$ from $0$ to $i$. We will use an inductive construction. At each point in the construction, we have a Toom cycle $(V,\vec E,\psi)$ rooted at $0$, and we let
\be
M:=\sup_{v\in V}{\rm dist}(\psi_v)
\ee
denote the largest distance from $0$ of all vertices of the Toom cycle. We will make sure that at each point in our construction, the following induction hypotheses are satisfied:
\begin{itemize}
\item[${\rm(i)'}$] if $\psi_v\not\in U_\ast$ for some $v\in V_\ast$, then ${\rm dist}(\psi_v)\geq M-1$,
\item[${\rm(ii)}$] $\dis(\psi_v,\psi_w)\in\vec H_s$ for all $(v,w)\in\vec E_s$ with $v\in V'_s$ $(s=1,2)$,
\item[${\rm(iii)}$] $\dis(\psi_v,\psi_w)\in\vec H_{3-s}$ for all $(v,w)\in\vec E_s$ with $v\in V'_\circ$ $(s=1,2)$.
\end{itemize}
Our construction will end as soon as in place of ${\rm(i)'}$ we have the stronger condition
\begin{itemize}
\item[{\rm(i)}] $\dis\psi_v\in U_\ast$ for all $v\in V_\ast$,
\end{itemize}
since this then guarantees that $(V,\vec E,\psi)$ is present in $(\La,\Hi)$. We note that in order to specify the Toom cycle, it suffices to know the function $\psi:[n]\to U$ only, since the induction hypotheses (ii) and (iii) imply that for $1\leq v\leq n$,
\be\label{etaupdo}
(v-1,v)\in\vec E\mbox{ if }(\psi_{v-1},\psi_v)\in\vec H
\quand
(v,v-1)\in\vec E\mbox{ if }(\psi_v,\psi_{v-1})\in\vec H.
\ee
Thus, at each step in the induction, we only specify an integer $n\geq 1$ and a function $\psi:[n]\to U$; it is then implicit that $\vec E$ is defined by (\ref{etaupdo}). It will be useful to view $\psi$ as a word $\psi_0\cdots\psi_n$ of length $n+1$, made up from the alphabet $U$, with $\psi_0=0=\psi_n$. We start with $n=1$ and $\psi_0=\psi_1:=0$. If $0\in U_\ast$ (and hence $|U|=1$), then we are done. In the opposite case, as long as (i) is not yet satisfied, we modify $\psi$ according to the following two steps, that are illustrated in Figure~\ref{fig:looperas}.

\begin{itemize}
\item[{\rm I.}] \emph{Exploration.} We pick a $v\in V_\ast$ such that $i:=\psi_v\not\in U_\ast$. If it is possible to pick $v$ such that ${\rm dist}(\psi_v)=M-1$, then we do so; in the opposite case we pick $v$ such that ${\rm dist}(\psi_v)=M$. By Proposition~\ref{P:typexpl}~(iii) there are unique $j_1,j_2\in U$ such that $(i,j_s)\in\vec G_s$ $(s=1,2)$. We modify the word $\psi_0\cdots\psi_n$ by inserting in place of the letter $\psi_v=i$ the string $i\,j_1\,i\,j_2\,i$.
\item[{\rm II.}] \emph{Loop erasion.} If as the result of the exploration, there are $v_1,v_2\in V_\ast$ with $v_1<v_2$ such that $i:=\psi_{v_1}=\psi_{v_2}$, then in place of the string $\psi_{v_1}\cdots\psi_{v_2}$ we insert the letter $i$. We repeat this until there are no more $v_1,v_2\in V_\ast$ with $v_1<v_2$ such that $\psi_{v_1}=\psi_{v_2}$.
\end{itemize}
We must check that at the end of each induction step, we obtain a Toom cycle satisfying the induction hypotheses ${\rm(i)'},{\rm(ii)},{\rm(iii)}$. We first investigate the effect of exploration.

After the exploration step, it is clear that $\psi$ via (\ref{etaupdo}) defines an oriented cycle. The map $\psi$ may no longer satisfy condition~(i) of Definition~\ref{D:cycle} (this will be fixed in the loop erasion step), but because of the way we have chosen $v$, after the exploration step, it will be true that:
\be\ba{l}\label{onlymax}
\mbox{if $\psi_w=\psi_{w'}$ for some $w\in V_\ast$ and $w'\in V$ with $w\neq w'$,}\\
\mbox{then ${\rm dist}(w)={\rm dist}(w')=M$ and $w'\in V_\ast$.} 
\ec
Indeed, the fact that ${\rm dist}(w)={\rm dist}(w')=M$ follows from the fact that the newly added vertices are at the largest distance $M$ from 0, while (\ref{VVV}) and (\ref{etaupdo}) imply that vertices at distance $M$ from $0$ must be elements of $V_\ast$. We claim that after the exploration step $\psi$ still satisfies condition (ii) of Definition~\ref{D:cycle}. Indeed, before the exploration step, (i) was still satisfied so the sink at $i$ did not overlap with any other vertices. In the exploration step, we add two sinks at the positions $j_1$ and $j_2$ and replace the old sink at the position $i$ by three vertices in $V_1,V_\circ$, and $V_2$, respectively, in this order. From this we see that after the exploration step, condition (ii) of Definition~\ref{D:cycle} is still satisfied.

We claim that after the exploration step, the induction hypotheses ${\rm(i)'},{\rm(ii)}$, and ${\rm(iii)}$ remain valid. Indeed, ${\rm(i)'}$ remains valid since $M$ does not increase unless all $v\in V_\ast$ for which $\psi_v\not\in U_\ast$ are at distance $M$ from $0$, and hence at least at distance $M-1$ after $M$ has increased. For the remaining induction hypotheses, we observe that in the exploration step, all existing edges of the oriented cycle keep their orientation. Their starting and endvertices also stay in whichever of the sets $V_\circ,V_\ast,V_1$, and $V_2$ they were in before, except (in the case $v\neq 0$) for the edges that ended in the vertex $v\in V_\ast$, whose new endvertices now belong to the sets $V_1$ and $V_2$, respectively. This has no influence on the induction hypotheses ${\rm(ii)}$ and ${\rm(iii)}$, however, for which only the starting vertices matter. Also, it is straightforward to check that the new edges inserted in the exploration step satisfy ${\rm(ii)}$ and ${\rm(iii)}$.

We next investigate the effect of loop erasion. During a loop erasion, all edges keep their charge and (because we are assuming $v_1,v_2\in V_\ast$) also all vertices stay in whichever of the sets $V_\circ,V_\ast,V_1$, and $V_2$ they were in before. Since they moreover preserve their relative order in $V$, this implies that condition (ii) of Definition~\ref{D:cycle} and the induction hypotheses ${\rm(ii)}$ and ${\rm(iii)}$ remain valid. In view of (\ref{onlymax}), during loop erasion, $M$ does not change and hence the induction hypothesis ${\rm(i)'}$ also remains valid. Furthermore, the process of loop erasion also restores condition~(i) of Definition~\ref{D:cycle}. This completes the induction step.

To complete the proof, we observe that in each step, either $M$ increases, or the number of vertices $v\in V_\ast$ with ${\rm dist}(\psi_v)=M-1$ and $\psi_v\not\in U_\ast$ decreases. By the finiteness of $(U,\Gi)$, this implies that our inductive construction terminates after a finite number of steps. It follows from our induction hypotheses and the fact that $(U,\Gi)$ is a subgraph of the typed dependence graph $(\La,\Hi)$ that at the end we obtain a Toom cycle that is present in $(\La,\Hi)$.
\epro

\section{The Peierls argument}\label{S:Peierls}

\subsection{Set-up}\label{S:setup}

In this section we prove Theorem~\ref{T:detspeed}, which gives sufficient conditions for the stability of monotone cellular automata with intrinsic randomness, as well as Proposition~\ref{P:coopbd}, which gives lower bounds on the critical noise parameter for two deterministic monotone cellular automata. Both proofs are based on the Peierls bound from Theorem~\ref{T:ToomPei}. In the present subsection we translate the setting of Theorem~\ref{T:detspeed} and Proposition~\ref{P:coopbd} into the more general language of Theorem~\ref{T:ToomPei} and make a choice for the typed dependence graph $(\La,\Hi)$ of Theorem~\ref{T:ToomPei} based on a linear polar function.

Throughout this section we assume that:
\begin{itemize}
\item $\phi_1,\ldots,\phi_m:\{0,1\}^{\Z^d}\to\{0,1\}$ are non-constant monotone local functions,
\item $\rbf=\big(\rbf(1),\ldots,\rbf(m)\big)$ is a probability distribution on $\{1,\ldots,m\}$.
\end{itemize}
For each $p\in[0,1]$, we let $\Phi^p=(\Phi^p_{i,t})_{(i,t)\in\Z^{d+1}}$ be an i.i.d.\ collection of maps as in (\ref{Phip}). We set $\La:=\Z^{d+1}$ and for each $(i,t)\in\La$ (with $i\in\Z^d$ and $t\in\Z$) we define $\Phi^p_{(i,t)}$ as in (\ref{Phip2}), so that $(\Phi^p_{(i,t)})_{(i,t)\in\La}$ is a random monotone  cellular automaton of the type considered in Theorem~\ref{T:ToomPei}. It will be convenient to define $\kappa:\La\to\{0,\ldots,m\}$ by
\be\label{kappa}
\Phi_{i,t}=:\phi_{\kappa(i,t)}\qquad\big((i,t)\in\La\big).
\ee
Then, in the notation of Theorem~\ref{T:ToomPei},
\be
\La_0=\big\{(i,t)\in\La:\kappa(i,t)=0\big\}
\quand
\La_\bullet=\big\{(i,t)\in\La:\kappa(i,t)\in\{1,\ldots,m\}\big\}.
\ee

In order to apply Theorem~\ref{T:ToomPei}, we need to choose ${\rm A}_{s,(i,t)}\in\Oi(\Phi_{(i,t)})$ for each $\lis$ and $(i,t)\in\La_\bullet$. We will let our choice be guided by a polar function. Throughout this section we assume that $L:\R^d\to\R^\sig$ is a linear polar function of dimension $\sig\geq 2$ such that
\be\label{epseps}
\eps:=\sum_{s=1}^\sig\eps_s>0\quad\mbox{with}\quad\eps_s:=\inf_{1\leq k\leq m}\eps_{\phi_k}(L_s)\quad(\lis),
\ee
where $\eps_{\phi_k}(L_s)$ is the edge speed defined in (\ref{edge}). For each $\lis$ and $1\leq k\leq m$, we fix $A_{s,k}\in\Oi(\phi_k)$ such that
\be\label{Ask}
\sup_{A\in\Oi(\phi_k)}\inf_{i\in A}L_s(i)=:\eps_{\phi_k}(L_s)=\inf_{i\in A_{s,k}}L_s(i) \quad (\lis,\ 1\leq k\leq m),
\ee
i.e., $A_{s,k}$ is a set for which the supremum in the definition of the edge speed in \eqref{edge} is achieved. Then setting
\be\label{Asit}
{\rm A}_{s,(i,t)}:=\big\{(i+j,t-1):j\in A_{s,\kappa(i,t)}\big\}\qquad\big(\lis,\ (i,t)\in\La_\bullet\big)
\ee
defines sets ${\rm A}_{s,(i,t)}\in\Oi(\Phi_{(i,t)})$ as needed for the application of Theorem~\ref{T:ToomPei}. In Subsection~\ref{S:edgebnd} below we will explain why this is the ``right'' choice for these sets. We let $(\La,\Hi)$ denote the typed dependence graph defined in terms of the sets ${\rm A}_{s,(i,t)}$ as in (\ref{Hs}), and $\big(\ov X^p(i,t)\big)_{(i,t)\in\La}$ denote the maximal trajectory of $\Phi^p$. We denote the origin in $\Z^{d+1}=\Z^d\times\Z$ by $(0,0)$, and write $\Ti_{(0,0)}$ for the set of Toom contours rooted at $(0,0)$, up to equivalence. Then Theorem~\ref{T:ToomPei} tells us that
\be\label{Pei}
1-\ov\rho(p)=\P\big[\ov X^p(0,0)=0\big]\leq\sum_{T\in\Ti_{(0,0)}}\P\big[T\mbox{ is present in }(\La,\Hi)\big].
\ee
We observe that as a consequence of properties (ii) and (iii) of Definition~\ref{D:present}, each Toom contour $T=(v_\circ,\Vi,\Ei,\psi)$ with $\sig$ charges that is present in $(\La,\Hi)$ must satisfy:
\be\ba{rl}\label{goodcontour}
\dis{\rm(ii)'}&\dis\psi(w)=\psi(v)+(j,-1)\mbox{ for some }j\in\De_s\mbox{ for all }(v,w)\in\vec E^\bullet_s\quad(\lis),\\[5pt]
\dis{\rm(iii)'}&\dis\psi(w)=\psi(v)+(j,-1)\mbox{ for some }j\in\De\mbox{ for all }(v,w)\in\vec E^\circ,
\ec
where
\be\label{DeDef}
\De_s:=\bigcup_{k=1}^mA_{s,k}\quad(\lis)\quand\De:=\bigcup_{s=1}^\sig\De_s.
\ee
We let $\Ti'_{(0,0)}$ denote the set of all $T\in\Ti_{(0,0)}$ that satisfy (\ref{goodcontour}). Then clearly, in (\ref{Pei}) we can restrict the sum to $T\in\Ti'_{(0,0)}$ since all other terms are zero.

For Toom cycles, similar arguments apply. In this case, we won't need the concept of equivalence of Toom contours defined in Definition~\ref{D:contour} but can use the slightly weaker but more intuitive concept of isomorphism of Toom contours. In line with this, we let $\bar\Ti_{(0,0)}$ denote the set of Toom cycles rooted at $(0,0)$, up to isomorphism, and inspired by Definition~\ref{D:cycpres}, we let $\bar\Ti'_{(0,0)}$ denote the subset of Toom cycles that moreover satisfy, for $s=1,2$,
\be\ba{rl}\label{goodcycle}
\dis{\rm(ii)''}&\dis\psi(w)=\psi(v)+(j,-1)\mbox{ for some }j\in\De_s\mbox{ for all }(v,w)\in\vec E_s\mbox{ with }v\in V'_s,\\[5pt]
\dis{\rm(iii)''}&\dis\psi(w)=\psi(v)+(j,-1)\mbox{ for some }j\in\De_{3-s}\mbox{ for all }(v,w)\in\vec E_s\mbox{ with }v\in V'_\circ.
\ec
Then Theorem~\ref{T:ToomPei} tells us that
\be\label{Peic}
1-\ov\rho(p)=\P\big[\ov X^p(0,0)=0\big]\leq\sum_{T\in\bar\Ti'_{(0,0)}}\P\big[T\mbox{ is present in }(\La,\Hi)\big].
\ee

\subsection{Stability of cellular automata with intrinsic randomness}\label{S:mainproof}

In this subsection we prove Theorem~\ref{T:detspeed}. For each Toom contour $T=(v_\circ,\mathcal V,\Ei,\psi)$ rooted at $(0,0)$ let
\be\label{nast}
n_\ast(T):=|V_\circ|=|V_\ast|
\ee
denote its number of sinks and sources, each. Recall that $\Ti'_{(0,0)}$ denotes the set of Toom contours with the additional properties (ii)' and (iii)' from (\ref{goodcontour}). The following lemma states that each $T\in\Ti'_{(0,0)}$ has an equal number of charged edges of each charge.

\bl[Number of charged edges]\label{L:edges}
For each Toom contour $T=(v_\circ,\Vi,\Ei,\psi)\in\Ti'_{(0,0)}$ with $\sig\geq 2$ charges there exists an integer $n_{\rm e}(T)$ such that
\be\label{nedges}
n_{\rm e}(T):=|\vec E_1|=\cdots=|\vec E_\sig|.
\ee
\el

\bpro
We write $\psi(v)=\big(\psi_1(v),\ldots,\psi_{d+1}(v)\big)$ where $\psi_{d+1}(v)$ denotes the time coordinate. The conditions in (\ref{goodcontour}) imply that $\psi_{d+1}(v)-\psi_{d+1}(w)=1$ for each $(v,w)\in\vec E$. Recall that by Definition~\ref{D:toomgraph} in a Toom graph at each source there emerge $\sig$ charges, one of each type, that then travel via internal vertices of the corresponding charge through the graph until they arrive at a sink, in such a way that at each sink there converge precisely $\sig$ charges. This implies
\be
|\vec E_1|=\cdots=|\vec E_\sig|=\sum_{v\in V_\ast}\psi_{d+1}(v)-\sum_{v\in V_\circ}\psi_{d+1}(v).
\ee
\epro

To prove Theorem~\ref{T:detspeed} we need two more lemmas, the proof of which will be postponed till later. To state the first lemma, let
\be\label{eq:Nn}
N_n:=\big|\{T\in\Ti'_{(0,0)}:n_{\rm e}(T)=n\}\big|\qquad(n\geq 0)
\ee
denote the number of non-equivalent contours in $\Ti'_{(0,0)}$ that have $n$ edges of each charge. In Subsection~\ref{S:expbd} we will prove the following exponential bound on $N_n$.

\bl[Exponential bound]
Let\label{L:expbd} $M:=\big|\De\big|$ with $\De$ defined in (\ref{DeDef}) and let $\tau:=\lceil\ha\sig\rceil$ denote $\ha\sig$ rounded up to the next integer. Then
\be
N_n\leq n^{\tau-1}(\tau+1)^{2\tau n}M^{\sig n}\qquad(n\geq 1).
\ee
\el

For our next lemma, we fix a polar function $L$ satisfying the assumptions of Theorem~\ref{T:detspeed} and we define
\be\ba{r@{\,}c@{\,}lcr@{\,}c@{\,}ll}\label{Rdef}
\dis R&:=&\dis\sum_{s=1}^\sig R_s&\quad\mbox{with}\quad&
\dis R_s&:=&\dis-\inf_{i\in\De}L_s(i)\qquad&\dis(1\leq s\leq\sig),
\ec
and we recall that $\eps$ and $\eps_s$ are defined in (\ref{epseps}). We will prove the following lemma in Subsection~\ref{S:edgebnd}.

\bl[Upper bound on the number of edges]
Each\label{L:edgebnd} Toom contour $T\in\Ti'_{(0, 0)}$ satisfies $n_{\rm e}(T)\leq(1+R/\eps)\big(n_\ast(T)-1\big)$.
\el

\bpro[of Theorem~\ref{T:detspeed}]
We use \eqref{Pei} which follows from Theorem~\ref{T:ToomPei}. To prove the stability of $\Phi^0$, it is enough to prove that the right-hand-side of \eqref{Pei} goes to 0 as $p\to 0$, while by the remarks below \eqref{Pei} it suffices to sum over all $T\in\Ti'_{(0,0)}$. By condition~(i) of Definition~\ref{D:embed} of an embedding, sinks of a Toom contour do not overlap. By condition~(i) of Definition~\ref{D:present} of what it means for a Toom contour to be present, each sink corresponds to a space-time point $(i,t)$ that is defective, meaning that $\Phi^p_{i,t}=\phi^0$, which happens with probability $p$, independently for all space-time points. As a result, the probability that a given contour $T$ is present in $(\La,\Hi)$ can be estimated from above by $p^{n_\ast(T)}$. By Lemma~\ref{L:edgebnd}, it follows that
\be\ba{l}\label{Peierls}
\dis 1-\ov\rho(p)\leq \sum_{T\in\Ti'_{(0,0)}}\P\big[T\mbox{ is present in }(\La, \mathcal H)\big]
\leq\sum_{T\in\Ti'_{(0,0)}}p^{n_\ast(T)}=p\sum_{T\in\Ti'_{(0,0)}}p^{n_\ast(T)-1}\\[5pt]
\dis\qquad\qquad\leq p\sum_{T\in\Ti'_{(0,0)}}p^{n_{\rm e}(T)/(1+R/\eps)}
=p\sum_{n=0}^\infty N_n p^{n/(1+R/\eps)},
\ec
Combining (\ref{Peierls}) and Lemma~\ref{L:expbd}, we see that this sum is finite for $p$ sufficiently small and hence (by dominated convergence) tends to zero as $p\to 0$. This proves that $\ov\rho(p)\to 1$ as $p\to 0$.
\epro

\subsection{Bounding the edges in terms of the sinks}\label{S:edgebnd}

In this subsection, we prove Lemma~\ref{L:edgebnd}, which says that for Toom contours in $\Ti'_{(0,0)}$, the number of edges can be bounded in terms of the number of sinks. Before we give the formal proof, we explain the main idea, which is really the central idea behind the proof of Theorem~\ref{T:detspeed} and the definition of Toom contours.

As explained in the previous subsection, the probability that a contour $T$ is present can be estimated from above by $p^{n_\ast(T)}$, where $n_\ast(T)$ is the number of sinks of the Toom contour. Therefore, we can estimate the expected number of Toom contours that is present in $(\La,\Hi)$ from above by $\sum_nM_np^n$, where $M_n$ denotes the number of non-equivalent contours in $\Ti'_{(0,0)}$ with $n$ sinks. In general, it is difficult to control the number of contours with a given number of sinks. As shown in Lemma~\ref{L:expbd}, however, we have good control over the number of contours with a given number of edges. Therefore, as we have seen in the proof of Theorem~\ref{T:detspeed}, to show that the Peierls sum in (\ref{Pei}) is small if $p$ is small, it suffices to have a result like Lemma~\ref{L:edgebnd} that bounds the number of edges from above in terms of the number of sinks.

It is precisely here that condition \eqref{eps} of Theorem~\ref{T:detspeed} on the worst-case edge speeds is used. Recall that $\vec E_s$ are the directed edges of charge $s$, which are distinguished as in (\ref{Ecirc}) into those that come out of a source other that the root (the set $\vec E^\circ_s$) and the others (the set $\vec E^\bullet_s$). Condition~(ii) of Definition~\ref{D:present} says that edges in $\vec E^\bullet_s$ must be embedded at edges of the same charge of the typed dependence graph $(\La,\Hi)$. Here $\vec H_s$ is defined in (\ref{Hs}) where the sets ${\rm A}_{s,i}$ are chosen in relation to the polar function $L$ as in (\ref{Ask}) and (\ref{Asit}). The result of all this is that:
\begin{itemize}
\item The function $L_s$ must increase by at least $\eps_s$ along each edge of charge $s$, except for edges that come out of sources other than the root.
\end{itemize}
Using this, condition \eqref{eps}, Lemma~\ref{L:edges}, and the fact that one edge of each charge originates at each source and one edge of each charge arrives at each sink, we can bound the number of edges in $\vec E^\bullet_s$ in terms of the number of edges in $\vec E^\circ_s$. Since there are $|\vec E^\circ_s|+1$ sources and an equal number of sinks, this allows us to bound the total number of edges in terms of the number of sinks.

We now make these ideas precise and prove Lemma~\ref{L:edgebnd}. We start with a general observation. On any set $\La$, we define a \emph{polar function} of \emph{dimension} $\sig\geq 2$ to be a function ${\rm L}:\La\to\R^\sig$ such that
\be\label{polar2}
\sum_{s=1}^\sig{\rm L}_s(i)=0\qquad(i\in\La).
\ee
The following lemma makes a connection between Toom contours and polar functions.

\bl[Zero sum property]
Let\label{L:zerosum} $(v_\circ,\Vi,\Ei,\psi)$ be a Toom contour with $\sig\geq 2$ charges and let ${\rm L}:\La\to\R^\sig$ be a polar function of dimension $\sig$. Then
\be\label{zerosum}
\sum_{s=1}^\sig\sum_{(v,w)\in\vec E_s}\big({\rm L}_s(\psi(w))-{\rm L}_s(\psi(v))\big)=0.
\ee
\el

\bpro
We can rewrite the sum in (\ref{zerosum}) as
\be
\sum_{v\in V}\Big\{\sum_{s=1}^\sig\sum_{(v,w)\in\vec E_{s,{\rm out}}(v)}{\rm L}_s(\psi(v))
-\sum_{s=1}^\sig\sum_{(u,v)\in\vec E_{s,{\rm in}}(v)}{\rm L}_s(\psi(v))\Big\}.
\ee
At internal vertices, the term inside the brackets is zero because the number of incoming edges of each charge equals the number of outgoing edges of that charge. At the sources and sinks, the term inside the brackets is zero by the defining property (\ref{polar2}) of a polar function, since there is precisely one outgoing (resp.\ incoming) edge of each charge.
\epro

\bpro[of Lemma~\ref{L:edgebnd}]
We trivially ``lift'' the linear polar function $L$, which is defined on $\Z^d$, to the space-time set $\La=\Z^{d+1}$ by setting
\be
{\rm L}_s(i,t):=L_s(i)\qquad\big(i\in\Z^d,\ t\in\Z).
\ee
Let $(v_\circ,\Vi,\Ei,\psi)=T\in\Ti'_{(0,0)}$. We claim that
\be\ba{ll}\label{Lincr}
\dis{\rm L}_s\big(\psi(w)\big)-{\rm L}_s\big(\psi(v)\big)\geq\eps_s
\quad&\dis\mbox{if }(v,w)\in\vec E^\bullet_s,\\[5pt]
\dis{\rm L}_s\big(\psi(w)\big)-{\rm L}_s\big(\psi(v)\big)\geq -R_s
\quad&\dis\mbox{if }(v,w)\in\vec E^\circ_s.
\ec
Indeed, by condition (ii)' in the definition of the set $\Ti'_{(0,0)}$ in \eqref{goodcontour}, $(v,w)\in\vec E^\bullet_s$ implies $\psi(w)=\psi(v)+(j,-1)$ for some $j\in\De_s=\bigcup_{k=1}^mA_{s,k}$. The linearity of $L_s$ implies that ${\rm L}_s(\psi(w))-{\rm L}_s(\psi(v))=L_s(j)$, which is $\geq\eps_s$ for all $j\in\De_s$ by (\ref{epseps}) and (\ref{Ask}). The second inequality in (\ref{Lincr}) follows in the same way from condition (iii)' in (\ref{goodcontour}) and (\ref{Rdef}).

By their definition in \eqref{Ecirc} and Lemma~\ref{L:edges}, we have 
\be\label{siso}
|\vec E^\circ_s|=n_\ast(T)-1 \quad \mbox{ and } \quad |\vec E^\bullet_s|=n_{\rm e}(T)-n_\ast(T)+1 \quad (1\leq s\leq\sig).
\ee
Lemma~\ref{L:zerosum}, (\ref{Lincr}), and (\ref{siso}) now imply that
\bc\label{siso2}
\dis 0
&=&\dis\sum_{s=1}^\sig\Big(\sum_{(v,w)\in \vec E^\bullet_s}
\big({\rm L}_s(\psi(w))-{\rm L}_s(\psi(v))\big)
+\sum_{(v,w)\in\vec E^\circ_s}
\big({\rm L}_s(\psi(w))-{\rm L}_s(\psi(v))\big)\Big)\\[5pt]
&\geq&\dis\sum_{s=1}^\sig\big[\big(n_{\rm e}(T)-n_\ast(T)+1\big)\eps_s-\big(n_\ast(T)-1\big)R_s\big]
=\eps n_{\rm e}(T)-(\eps+R)\big(n_\ast(T)-1\big),
\ec
which implies $n_{\rm e}(T)\leq(1+R/\eps)\big(n_\ast(T)-1\big)$.
\epro

\subsection{Exponential bounds on the number of contours}\label{S:expbd}

In this subsection, we provide the only missing ingredient in the proof of Theorem~\ref{T:detspeed}, which is the proof of Lemma~\ref{L:expbd}. If we would be satisfied with just any exponential bound, then the proof could be quite short, but with a view towards Proposition~\ref{P:coopbd} we will argue a bit more carefully to get a sharper bound. The idea of the proof is to walk around in a Toom contour in such a way that each edge is traversed at least once, and to record enough information along the way to be able to uniquely reconstruct the Toom contour.\med

\bpro[of Lemma~\ref{L:expbd}]
We first consider the case that the number of charges $\sig$ is even. Let $T=(v_\circ,\Vi,\Ei,\psi)\in\Ti'_{(0,0)}$. Recall that $(\Vi,\Ei)$ is a typed directed graph with $\sig$ types of edges, that are called charges. In $(\Vi,\Ei)$, all edges point in the direction from the sources to the sinks. We modify $(\Vi,\Ei)$ by reversing the direction of edges of the charges $\ha\sig+1,\ldots,\sig$. Let $(\Vi,\Ei')$ denote the modified graph. In $(\Vi,\Ei')$, the number of incoming edges at each vertex equals the number of outgoing edges. Since moreover the undirected graph $(V,E)$ is connected, it is not hard to see\footnote{This is a simple variation of the ``Bridges of K\"onigsberg'' problem that was solved by Euler.} that it is possible to walk through the directed graph $(\Vi,\Ei')$ starting from the root using an edge of charge $1$, in such a way that each directed edge of $\Ei'$ is traversed exactly once.

Let $l:=\sig n_{\rm e}(T)$ denote the total number of edges of $(\Vi,\Ei')$ and for $0<k\leq l$, let $(v_{k-1},v_k)\in\vec E'_{s_k}$ denote the $k$-th step of the walk, which has charge $s_k$. Write $\psi(v_k)=:\psi(v_{k-1})+(\de_k,\pm 1)$ where $\de_k$ is the spatial increment of the $k$-th step and $\pm 1$ is the temporal increment, which is determined by the charge $s_k$ of the $k$-th step: it is -1 for charges $1, \dots, \frac 1 2 \sig$ and 1 otherwise. Let $k_0,\ldots,k_{\sig/2}$ denote the times when the walk visits the root $v_\circ$. We claim that in order to specify $(v_\circ,\Vi,\Ei,\psi)$ uniquely up to equivalence, in the sense defined in Definition~\ref{D:contour}, it suffices to know the sequences
\be
(s_1,\ldots,s_l),\quad(\de_1,\ldots,\de_l),\quand(k_0,\ldots,k_{\sig/2}).
\ee
Indeed, the sinks and sources correspond to changes in the temporal direction of the walk which can be read off from the charges. Although the images under $\psi$ of sources may overlap, we can identify which edges connect to the root, and we also know the increment of $\psi(v_k)-\psi(v_{k-1})$ in each step, hence all objects in (\ref{psiEi}) can be identified.

The first charge $s_1$ is 1 and after that, in each step, we have the choice to either continue with the same charge or choose one of the other $\ha\sig$ available charges. This means that there are no more than $(\ha\sig+1)^{l-1}$ possible ways to specify the charges $(s_1,\ldots,s_l)$. Recalling $M=|\Delta|=\big|\bigcup_{s=1}^\sig \bigcup_{k=1}^m A_{s,k}\big|$, we see that there are no more than $M^l$ possible ways to specify the spatial increments $(\de_1,\ldots,\de_l)$. Since $k_0=0,k_{\sig/2}=l$, we can roughly estimate the number of ways to specify the visits to the root from above by $n^{\sig/2-1}$. Recalling that $l=\sig n_{\rm e}(T)$, this yields the bound
\be
N_n\leq n^{\sig/2-1}(\ha\sig+1)^{\sig n-1}M^{\sig n}.
\ee
This completes the proof when $\sig$ is even.

When $\sig$ is odd, we modify $(\Vi,\Ei)$ by doubling all edges of charge $\sig$, i.e., we define $(\Vi,\Fi)$ with
\be
\vec F=(\vec F_1,\ldots,\vec F_{\sig+1}):=(\vec E_1,\ldots,\vec E_\sig,\vec E_\sig),
\ee
and next we modify $(\Vi,\Fi)$ by reversing the direction of all edges of the charges $\lceil\ha\sig\rceil+1,\ldots,\sig+1$. We can define a walk in the resulting graph $(\Vi,\Fi')$ as before and record the charges and spatial increments for each step, as well as the visits to the root. In fact, in order to specify $(v_\circ,\Vi,\Ei,\psi)$ uniquely up to equivalence, we do not have to distinguish the charges $\sig$ and $\sig+1$. Recall that edges of the charges $\sig$ and $\sig+1$ result from doubling the edges of charge $\sig$ and hence always come in pairs, connecting the same vertices. Since sinks do not overlap and internal vertices of a given charge do not overlap, and since we traverse edges of the charges $\sig$ and $\sig+1$ in the direction from the sinks towards the sources, whenever we are about to traverse an edge that belongs to a pair of edges of the charges $\sig$ and $\sig+1$, we know whether we have already traversed the other edge of the pair. In view of this, for each pair, we only have to specify the spatial displacement at the first time that we traverse an edge of the pair. Using these considerations, we arrive at the bound
\be
N_n\leq n^{\lceil\sig/2\rceil-1}(\lceil\ha\sig\rceil+1)^{(\sig+1)n-1}M^{\sig n}.
\ee
\epro

\subsection{Some bounds for Toom cycles}\label{S:twochar}

With Theorem~\ref{T:detspeed} proved, we start to prepare for the proof of Proposition~\ref{P:coopbd}. In the present subsection, we prove more precise versions of Lemmas \ref{L:expbd} and \ref{L:edgebnd} that hold only for Toom cycles and that will help us to get a better bound for the critical noise parameter $p_{\rm c}$ of the cellular automaton defined by the map $\phi^{\rm coop}$. Recall that at the end of Subsection~\ref{S:setup} we denoted the set of non-isomorphic Toom cycles rooted at $(0,0)$ by $\bar\Ti_{(0,0)}$ and we wrote $\bar\Ti'_{(0,0)}$ for the set of $T\in\bar\Ti_{(0,0)}$ that satisfy (\ref{goodcycle}). 

Similarly to (\ref{eq:Nn}), we let
\be\label{barN}
\bar N_n:=\big|\{T\in\bar\Ti'_{(0,0)}:n_{\rm e}(T)=n\}\big|\qquad(n\geq 0)
\ee
denote the number of non-isomorphic Toom cycles in $\bar\Ti'_{(0,0)}$ that have $n$ edges of each charge. We then have the following analogue of Lemma~\ref{L:expbd}.

\bl[Exponential bound for $\sig=2$]
Let\label{L:expbdcycle} $M_s:=|\De_s|$ $(s=1,2)$ with $\De_s$ defined in (\ref{DeDef}). Then
\be
\bar N_n\leq\ha(4M_1M_2)^{n} \qquad(n\geq 1).
\ee
\el

\bpro
The proof goes along the same lines as that of Lemma~\ref{L:expbd} for the case $\sig$ is even. Observe that for $\sig=2$, the walk visits the root 0 twice: $k_0=0, k_1=l$, where $l$ is the total number of edges of the cycle. Thus $(k_0, k_1)$ is deterministic, and we only need to specify the sequences
\be
(s_1,\ldots,s_l),\quad(\de_1,\ldots,\de_l).
\ee
Note that in this case, these sequences determine the Toom cycle up to isomorphism and not only up to equivalence as in the proof of Lemma~\ref{L:expbd}. The first charge $s_1$ is 1 and after that, in each step, we have the choice to either continue with the same charge or choose charge 2. This means that there are no more than $2^{l-1}$ possible ways to specify the charges $(s_1,\ldots,s_l)$. Once we have done that, by condition~(iii)'' of (\ref{goodcycle}), we know for each $0<k\leq l$ whether the spatial increment $\de_k$ is in $\De_1$ or $\De_2$. Recalling $M_s=|\De_s|$ $(s=1,2)$ and using the fact that $|\vec E_1|=|\vec E_2|=n_{\rm e}(T)=l/2$, we see that there are no more than $M_1^{l/2} \cdot M_2^{l/2}$ possible ways to specify $(\de_1,\ldots,\de_l)$. This yields the bound
\be
\bar N_n\leq 2^{2n-1}M_1^{n} \cdot M_2^{n}.
\ee
\epro

From now on, we fix a polar function $L$ of dimension $2$ satisfying the assumptions of Theorem~\ref{T:detspeed}. In analogy with (\ref{Rdef}), but with a view towards (\ref{goodcycle}) which in the present context replaces (\ref{goodcontour}), we define
\be\label{eq:Rprime}
\bar R:=\sum_{s=1}^2\bar R_s\qquad\mbox{with}\quad\bar R_1:=-\inf_{i\in\De_2}L_1(i)\quad\mbox{and}\quad\bar R_2:=-\inf_{i\in\De_1}L_2(i).
\ee
The following lemma is similar to Lemma~\ref{L:edgebnd}.

\bl[Upper bound on the number of edges for $\sig=2$]
Let\label{L:edgebndcycle} $\eps$ be defined in \eqref{eps} and let $\bar R$ be defined in \eqref{eq:Rprime}. Then each $T\in\bar\Ti'_{(0,0)}$ satisfies $n_{\rm e}(T)\leq(1+\bar R/\eps)\big(n_\ast(T)-1\big)$.
\el

\bpro
The proof is the same as that of Lemma~\ref{L:edgebnd}, with the only difference that condition (iii)'' of (\ref{goodcycle}) allows us to use $\bar R_s$ instead of $R_s$ $(s=1,2)$ as upper bounds.
\epro

\subsection{Finiteness of the Peierls sum}\label{S:finP}

We continue our preparations for the proof of Proposition~\ref{P:coopbd}. Our aim is to derive a lower bound $p_\ast$ on the critical noise parameter $p_{\rm c}$ defined in (\ref{pc}), which requires us to prove that $\ov\rho(p)>0$ for all $p<p_\ast$. By (\ref{Pei}), we have $\ov\rho(p)>0$ as soon as the Peierls sum on the right-hand side of (\ref{Pei}) is less than one. In the present subsection, we will prove that in fact it (more or less) suffices to show that the Peierls sum is finite. This will not only lead to slightly better bounds but also simplify our calculations later. Similar results, which say that finiteness of the Peierls sum already implies a phase transition, have been proved before. For percolation, the argument is quite simple \cite[Section~6a]{Dur88} but for other models such results can be a bit harder to obtain \cite{KSS14}.

We will work in the set-up of Subsection~\ref{S:setup}, but specialised to the case $m=1$, which means that $\Phi^0$ is a deterministic monotone cellular automaton. To simplify notation, we set $\phi:=\phi_1$ and $A_s:=A_{s,1}$. We recall from (\ref{Ask}) that $A_s$ is chosen such that
\be\label{Achoice}
\sup_{A\in\Oi(\phi)}\inf_{i\in A}L_s(i)=\inf_{i\in A_s}L_s(i)\qquad(\lis),
\ee
and that by (\ref{edge}) and (\ref{epseps})
\be\label{neweps}
\eps=\sum_{s=1}^\sig\eps_s\quad\mbox{with}\quad\eps_s=\eps_\phi(L_s)=\inf_{i\in A_s}L_s(i)\quad(\lis).
\ee
In our present setting, the definition of the sets $\De_s$ and $\De$ from (\ref{DeDef}) simplifies to
\be
\De_s=A_s\quad(\lis)\quand\De=\bigcup_{s=1}^\sig A_s.
\ee
As in Subsection~\ref{S:setup}, we let $\Ti'_{(0,0)}$ denote the set of Toom contours rooted at $(0,0)$ that satisfy (\ref{goodcontour}). In the special case $\sig=2$, we let $\bar\Ti'_{(0,0)}$ denote the set of Toom cycles rooted at $(0,0)$ that satisfy (\ref{goodcycle}). As in (\ref{nast}) we let $n_\ast(T)$ denote the number of sinks of a Toom contour $T$, which equals the number of sources. Here is the main result of this subsection.

\bp[Finiteness of the Peierls sum]
Assume\label{P:Peifin} that $\eps>0$. Then the condition
\be\label{Peifin}
\sum_{T\in\Ti'_{(0,0)}}p^{n_\ast(T)}<\infty
\ee
implies that $\ov\rho(p)>0$. If $\sig=2$, then the same conclusion can be drawn if in (\ref{Peifin}) the sum over $\Ti'_{(0,0)}$ is replaced by the sum over $\bar\Ti'_{(0,0)}$.
\ep

Proposition~\ref{P:Peifin} actually stops short of what we promised, since the sum in (\ref{Peifin}) is only an upper bound for the Peierls sum on the right-hand side of (\ref{Pei}). For our purposes, the statement of Proposition~\ref{P:Peifin} will be sufficient, however.

The proof of Proposition~\ref{P:Peifin} needs some preparations. Recall that $\La=\Z^{d+1}$ and that the function $\kappa$ is defined in (\ref{kappa}). In the present setting, $\big(\kappa(i,t)\big)_{(i,t)\in\La}$ are i.i.d.\ $\{0,1\}$-valued random variables with $\P[\kappa(i,t)=1]=p$, and
\be
\La_0:=\big\{(i,t)\in\La:\kappa(i,t)=0\big\}
\quand
\La_\bullet:=\big\{(i,t)\in\La:\kappa(i,t)=1\big\}.
\ee
Recall that sites in $\La_0$ are called \emph{defective}. In our present setting, the definition of the typed dependence graph $(\La,\Hi)$ simplifies to
\be\label{vecHs}
\vec H_s:=\big\{\big((i,t),(i+j,t-1)\big):(i,t)\in\La_\bullet,\ j\in A_s\big\}\qquad(\lis).
\ee
It will be convenient to define a modified typed dependence graph $(\La,\Hi^<)$ that has no defective sites $(i,t)$ with time coordinates $t>0$. Formally, we define $\La^<_0:=\big\{(i,t)\in\La_0:t\leq 0\big\}$, $\La^<_\bullet:=\La\beh\La^<_0$, and we define $\Hi^<=(\vec H^<_1,\ldots,\vec H^<_\sig)$ as in (\ref{vecHs}) but with $\La_\bullet$ replaced by $\La^<_\bullet$. We let $\Psi^<$ denote the monotone cellular automaton associated with $(\La,\Hi^<)$, in the sense of Definition~\ref{D:typdep}, and we let $\ov Y^<$ denote the maximal trajectory of $\Psi^<$.

\bl[Presence of a large contour]
Fix\label{L:largcont} $j_s\in A_s$ $(\lis)$ and $r\in\N$. Let $C_r\sub\Z^d$ with $r\in\N$ be inductively defined by
\be
C_0:=\{0\}\quand C_{r+1}:=\big\{i+j_s:i\in C_r,\ \lis\big\}\qquad(r\geq 0).
\ee
Then on the event that $\ov Y^<(i,0)=0$ for all $i\in C_r$, there is a Toom contour $(v_\circ,\Vi,\Ei,\psi)$ rooted at $(0,r)$ present in $(\La,\Hi^<)$. If $\sig=2$, then a Toom cycle rooted at $(0,r)$ is present in $(\La,\Hi^<)$.
\el

\bpro
Since there are no defective sites with positive time coordinates, by Definition~\ref{D:typdep} and (\ref{vecHs}) we have
\be
\ov Y^<(i,t)=\bigvee_{s=1}^\sig\bigwedge_{j\in A_s}Y^<(i+j,t-1)\qquad(i\in\Z^d,\ t>0).
\ee
Using this and the assumption that $\ov Y^<(i,0)=0$ for all $i\in C_r$ we see by induction that $\ov Y^<(i,t)=0$ for all $i\in C_{r-t}$ $(0\leq t\leq r)$ and hence in particular $\ov Y^<(0,r)=0$. The claim now follows from Theorems \ref{T:contour} and \ref{T:cycle}.
\epro

\bl[Many sinks]
Let\label{L:manysink} $\eps$ and $R$ be defined as in (\ref{eps}) and (\ref{Rdef}). Assume that $T$ is a Toom contour rooted at $(0,r)$ that is present in $(\La,\Hi^<)$. Then $n_\ast(T)\geq r\eps/R+1$.
\el

\bpro
This follows from an argument similar to the proof of Lemma~\ref{L:edgebnd}. Since $(\La,\Hi^<)$ has no defective sites with positive time coordinates, any Toom contour $T=(v_\circ,\Vi,\Ei,\psi)$ that is rooted at $(0,r)$ and present in $(\La,\Hi^<)$ must satisfy $|\vec E^\bullet_s|\geq r$ $(\lis)$, so a calculation as in (\ref{siso2}) gives
\be
0\geq\sum_{s=1}^\sig\big[r\eps_s-(n_\ast(T)-1)R_s\big]=r\eps-(n_\ast(T)-1)R.
\ee
\epro

\bpro[of Proposition~\ref{P:Peifin}]
Let $\Psi$ denote the monotone cellular automaton associated with $(\La,\Hi)$. Let $\ov X^p$ denote the maximal trajectory of $\Phi^p$ and let $\ov Y$ denote the maximal trajectory of $\Psi$. As in the proof of Theorem~\ref{T:ToomPei} we see that $\ov Y\leq \ov X^p$ (pointwise). Moreover, since $(\La,\Hi)$ and $(\La,\Hi^<)$ agree up to time zero, it is easy to see that $\ov Y(i,t)=\ov Y^<(i,t)$ for all $t\leq 0$. Let $M_n$ denote the number of non-equivalent contours in $\Ti'_{(0,0)}$ with $n$ sinks. Lemmas \ref{L:largcont} and \ref{L:manysink} tell us that
\be\label{allzero}
\P\big[\ov Y^<(i,0)=0\ \forall i\in C_r\big]\leq\sum_{n=r\eps/R+1}^\infty M_np^n,
\ee
and the same is true with $\ov Y^<$ replaced by $\ov Y$, since they are equal at time zero. Using (\ref{Peifin}) and the assumption that $\eps>0$, we see that we can choose $r$ large enough such that the probability in (\ref{allzero}) is less than one. It follows that
\be
0<\P\big[\ov Y(i,0)=1\mbox{ for some }i\in C_r\big]
\leq|C_r|\cdot\P\big[\ov Y(0,0)=1\big]\leq|C_r|\cdot\P\big[\ov X^p(0,0)=1\big],
\ee
proving that $\ov\rho(p)>0$. The proof for Toom cycles is the same.
\epro

\subsection{Bounds on the critical noise parameter}\label{S:explic}

In this subsection we prove Proposition~\ref{P:coopbd}. We work in the set-up of Subsection~\ref{S:setup}, specialised to the case $m=1$, as summarised in Subsection~\ref{S:finP}.\med

\bpro[of Proposition~\ref{P:coopbd}]
First we consider the cellular automaton that applies the map $\phi^{{\rm coop}}$ (recall \eqref{phiNEC}) at each space-time point. We chose $\sig:=2$, and the polar function
\be
L_1(z):=-z_1-z_2, \quad L_2(z):=z_1+z_2 \quad (z_1, z_2)\in\R^2.
\ee
Recalling the minimal one-sets of $\phi^{{\rm coop}}$ from \eqref{ANEC}, we then choose $A_1,A_2\in\Oi(\phi^{\rm coop})$ satisfying (\ref{Achoice}) by setting $A_1:=\{(0,0)\}$, and $A_2:=\{(0,1), (1,0)\}$. This has the result that the constants from \eqref{eps}, \eqref{Rdef} and \eqref{eq:Rprime} are given by $\eps=1$, $R=1$ and $\bar R=1$. We first give a bound using Toom contours. Applying Lemma~\ref{L:expbd} with $M=3$, $\sig=2$ and $\tau=1$, the Peierls bound (\ref{Peierls}) gives
\be\label{coopPei}
1-\ov\rho(p)\leq\sum_{T\in\Ti'_{(0,0)}}p^{n_\ast(T)}\leq p\sum_{n=0}^\infty N_np^{n/(1+R/\eps)}\leq
p\sum_{n=0}^\infty2^{2\tau n}3^{2n}p^{n/2}.
\ee
By Proposition~\ref{P:Peifin}, to prove that $\ov\rho(p)>0$, it actually suffices to prove that the right-hand side of (\ref{coopPei}) is finite, which happens when $36p^{1/2}<1$, leading to the bound $p_{\rm c}\geq 36^{-2}$.

Since $\sig=2$, we can use Toom cycles instead. Using \eqref{Pei}, Lemma~\ref{L:expbdcycle} with $M_1=1$, $M_2=2$, and Lemma~\ref{L:edgebndcycle}, we find that
\be\label{coopPei2}
1-\ov\rho(p)\leq\sum_{T\in\bar\Ti'_{(0,0)}}p^{n_\ast(T)}\leq p\sum_{n=0}^\infty\bar N_np^{n/(1+\bar R/\eps)}\leq p+\ha p\sum_{n=1}^\infty8^np^{n/2}.
\ee
Again, by Proposition~\ref{P:Peifin}, it suffices to prove that the right-hand side is finite, which happens when $8 p^{1/2}<1$, leading to the bound $p_{\rm c}\geq 1/64$.

Now consider the cellular automaton that applies the map $\phi^{{\rm NEC}}$ (recall \eqref{phiNEC}) at each space-time point. We chose $\sig:=3$, and the polar function
\be\label{eq:Toompolar2}
L_1(z_1,z_2):=-z_1,\quad
L_2(z_1,z_2):=-z_2,\quad
L_3(z_1,z_2):=z_1+z_2.
\ee
Recalling the minimal one-sets of $\phi^{{\rm NEC}}$ from \eqref{ANEC} we choose $A_1,A_2,A_3\in\Oi(\phi^{\rm NEC})$ satisfying (\ref{Achoice}) by setting $A_1:=\{(0,0),(0,1)\}$, $A_2:=\{(0,0),(1,0)\}$, and $A_3:=\{(0,1),(1,0)\}$. This has the result that the constants from~\eqref{eps} and \eqref{Rdef} are given by $\eps=1$ and $R=2$. Applying Lemma~\ref{L:expbd} with $M=3$, $\sig=3$ and $\tau=2$, the Peierls bound (\ref{Peierls}) gives
\be
1-\ov\rho(p)
\leq\sum_{T\in\Ti'_{(0,0)}}p^{n_\ast(T)}
\leq p\sum_{n=0}^\infty n3^{4n}3^{3n}p^{n/3}.
\ee
By Proposition~\ref{P:Peifin}, it suffices to prove that the right-hand side is finite, which happens when $3^7p^{1/3}<1$, leading to the bound $p_{\rm c}\geq 3^{-21}$.
\epro

\appendix

\section{Bootstrap percolation}\label{A:boot}

Our simplification of Toom's argument has successfully been used in the study of bootstrap percolation in \cite{HS22}. In this appendix, we briefly elaborate on this connection.

Instead of perturbing a deterministic monotone cellular automaton with noise that is i.i.d.\ in space and time, one can also look at perturbations that are i.i.d.\ in space but constant in time. More precisely, if $\phi:\{0,1\}^{\Z^d}\to\{0,1\}$ is a monotone local map that is not constant, then it is interesting to look at i.i.d.\ collections of random variables $\Phi^p=(\Phi^{p}_i)_{i\in\Z^d}$ with values in $\{\phi^0,\phi\}$ such that
\be
\P\big[\Phi^p_i=\phi^0\big]=p\quand\P\big[\Phi^p_i=\phi\big]=1-p,
\ee
and modify the evolution in (\ref{Markov}) in the sense that the same map $\Phi^p_i$ is applied at each time $t>0$. We let $\ti\rho_\phi(p)$ denote the long-time limit of the density started from $\un 1$ for this type of evolution. Recall from  (\ref{ovrho}) that $\ov\rho_\phi(p)$ denotes the limiting density when the noise is i.i.d.\ in space and time.

\begin{defi}[Forms of stability]
We\label{D:boot} say that $\phi$ is \emph{stable in the bootstrap sense} if $\lim_{p\to 0}\ti\rho_\phi(p)=1$ and \emph{stable in Toom's sense} if $\lim_{p\to 0}\ov\rho_\phi(p)=1$.
\end{defi}

Toom's stability theorem (Theorem~\ref{T:Toom}) completely answers the question which monotone local maps $\phi$ are stable in Toom's sense. The analogue question for stability in the bootstrap sense has been answered more recently. In order to formulate a precise result, we first need to translate the problem into the language of bootstrap percolation. In analogy with notation introduced in (\ref{Psiphi}), we let $\Psi^p_\phi:\{0,1\}^{\Z^d}\to\{0,1\}^{\Z^d}$ denote the random map defined as
\be\label{Psiphip}
\Psi^p_\phi(x)(i):=\Phi_i\big((x(i+j))_{j\in\Z^d}\big)\qquad\big(x\in\{0,1\}^{\Z^d}\big).
\ee
Since the noise is the same in each time step,
\be\label{bootdens}
\ti\rho_\phi(p):=\lim_{t\to\infty}\P\big[(\Psi^p_\phi)^t(\un 1)(i)=1\big]\qquad(i\in\Z^d),
\ee
where $(\Psi^p_\phi)^t$ denotes the $t$-th iterate of the map $\Psi^p_\phi$. Note that if $(\Psi^p_\phi)^t(\un 1)(i)=0$ for some $i\in\Z^d$ and $t\geq 1$, then $(\Psi^p_\phi)^{t+s}(\un 1)(i)=0$ for all $s\geq 0$.

Recall that $\Oi(\phh)$ and $\Zi(\phh)$ denote the sets of minimal one-sets and zero-sets of a monotone local map $\phh$, defined in Subsections \ref{S:intro} and \ref{S:minexpl}. For any monotone local map $\phi$ we define a map $\ov\phi:\{0,1\}^{\Z^d}\to\{0,1\}$ in terms of its set of minimal zero-sets as
\be\label{eq:ovphi}
\Zi(\ov\phi):=\big\{Z: Z\in\Zi(\phi), 0\notin Z \big\}.
\ee
We say that two monotone local maps $\phi_1$ and $\phi_2$ from $\{0,1\}^{\Z^d}$ to $\{0,1\}$ are \emph{equivalent in the bootstrap sense} if $\ov\phi_1=\ov\phi_2$. Observe that in this case $\phi_1(x)=\phi_2(x)$ for all $x\in\{0,1\}^{\Z^d}$ such that $x(0)=1$. It is easy to see that if $\phi_1$ and $\phi_2$ are equivalent in the bootstrap sense, then $(\Psi^p_{\phi_1})^t(\un 1)$ and $(\Psi^p_{\phi_2})^t(\un 1)$ are equal in law, and moreover almost surely equal if we couple both processes in the obvious way, by applying the zero map in the same space-time points. Therefore, in view of (\ref{bootdens}), for equivalent maps $\ti\rho_{\phi_1}(p)=\ti\rho_{\phi_2}(p)$ $(0\leq p\leq 1)$. In particular, for each map $\phi$
\be
\ti\rho_\phi(p)=\ti\rho_{\ov\phi}(p).
\ee
The set $\Zi(\ov\phi)$ is called the \emph{update family} in the bootstrap percolation literature, and it is traditionally in terms of this set that the dynamics are described.

Bootstrap percolation was first introduced in \cite{CLR79} to model magnetic materials at low temperature, and has since been extensively studied (see~\cite{Mor17} for a review). The first stability result on $\Z^d$ was established in \cite{Sch92} for the update family consisting of all sets containing exactly $r$ neighbours of the origin. Presently, there is a complete characterisation of bootstrap percolation maps on $\Z^d$. Based on the geometry of the sets in the update family, bootstrap maps $\phi$ fall into three universality classes: maps $\phi$ in the supercritical and critical classes are completely unstable in the bootstrap sense meaning that $\ti\rho_\phi(p)=0$ for all $p>0$, while those in the subcritical class are stable in the bootstrap sense. These results were first established in~\cite{BSU15, BBPS16} in two dimensions and recently extended to higher dimensions in a series of papers~\cite{BBMS22a, BBMS22b, BBMS24}. It was shown in \cite{HS22} that the stability result in the subcritical class can also (and more simply) be obtained using Toom contours, marking the first application of the results presented in this paper. Additionally, \cite{HT24} applies Toom contours to establish exponential decay results for the largest cluster of zeros in the final configuration of subcritical bootstrap percolation, as well as for the largest cluster of zeros in space-time for monotone cellular automata defined by an eroder.

The following theorem completely answers the question which monotone local maps $\phi$ are stable in the bootstrap sense.

\bt[Stability in the bootstrap sense]
Let\label{T:BPstablecrit} $\phi:\{0,1\}^{\Z^d}\to\{0,1\}$ be a monotone local map. Then $\phi$ is stable in the bootstrap percolation sense if and only if there exists a linear polar function $L$ of dimension $\sig\geq 2$ such that the edge speeds defined in (\ref{edge}) of the map $\ov\phi$ from \eqref{eq:ovphi} satisfy
\be\label{stableBP}
\eps_{\ov\phi}(L_s)>0\quad(\lis).
\ee
\et

This theorem is proved in~\cite[Lemma~2.1]{HS22}. (The lemma only states one direction of the theorem, but it is easy to check that every implication used in the proof is in fact an equivalence.) We observe that (\ref{stableBP}) is stronger than the condition that $\sum_{s=1}^\sig\eps_{\ov\phi}(L_s)>0$, so comparing with Theorem~\ref{T:Toom} and Lemma~\ref{L:erode} we see that stability in the bootstrap sense implies that $\ov\phi$ is stable in Toom's sense but not vice versa. To see a concrete example, consider the Duarte model on $\Z^2$ defined by the monotone local map
\be\label{phiDuarte}
\dis\phi^{\rm Duarte}(x):=\dis{\tt round}\big([x(0,1)+x(0,-1)+x(-1,0)]/3\big).
\ee
This rule is unstable in the bootstrap sense as it is known to belong to the critical universality class. However, as the intersection of the convex hulls of its minimal one-sets is empty, $\phi^{\rm Duarte}$ is an eroder and hence stable in Toom's sense. For any monotone local map $\phi$ we can apart from the map $\ov\phi$ from (\ref{eq:ovphi}) also define
\[
\Zi(\un\phi):=\Zi(\ov\phi)\cup\big\{ \{0\}\big\}.
\]
Then $\un\phi$ is equivalent to $\phi$ in the bootstrap sense\footnote{In fact, it is easy to see that $\un\phi$ and $\ov\phi$ are the smallest and largest maps that are equivalent to $\phi$.} while it is easy to see that $\un\phi$ is never stable in Toom's sense, as once a site flips to 0 it forever remains in state 0. Thus, for general monotone local maps $\phi$ all four combinations (stable/unstable in Toom's/the bootstrap sense) are possible.

\subsection*{Acknowledgement}

We thank Anja Sturm who was involved in the earlier phases of writing this paper for her contributions to the discussions. We thank Ivailo Hartarsky for useful discussions.

\subsection*{Declarations}

\emph{Conflict of interest statement} The first author is supported by grant 20-08468S of the Czech Science Foundation (GA\v{C}R). The second and third authors are supported by European Research Council Starting Grant 680275 ``MALIG''.\med

\noi
\emph{Data availability statement} The manuscript has no associated data.

\end{document}